\documentclass[reqno, 11pt]{amsart}
\usepackage[utf8]{inputenc}
\usepackage[T1]{fontenc}
\usepackage{lmodern}
\usepackage{epsfig,amsmath, amsthm ,amsfonts,latexsym,stmaryrd, amssymb}
\usepackage[abs]{overpic}
\usepackage[usenames,dvipsnames]{xcolor}
\usepackage[linktocpage=true]{hyperref}
\usepackage{caption}
\usepackage{subcaption}
\usepackage{dsfont}
\usepackage{mathtools}
\usepackage{tikz-cd}
\usepackage{simpler-wick}
\usetikzlibrary{shapes,arrows}
\usepackage{comment}
\usepackage{enumitem}
\usepackage{marginnote}
\usepackage[all,cmtip]{xy}
\usepackage[backend=biber, maxnames=4, style=alphabetic]{biblatex}
\hypersetup{
  colorlinks,
  citecolor=teal,
  linkcolor=purple,
  urlcolor=teal}
\theoremstyle{plain}

\newtheorem{theorem}{Theorem}

\newtheorem{dummy}{anything}[section]
\newtheorem{thm}{Theorem}

\newtheorem{proposition}{Proposition}[section]

\theoremstyle{definition}
\newtheorem{definition}[dummy]{Definition}

\newtheorem{example}[dummy]{Example}

\newtheorem{remark}[dummy]{Remark}

\theoremstyle{remark}
\newcommand{\M}{\mathcal{M}}

\DeclareFontFamily{U}{mathx}{}
\DeclareFontShape{U}{mathx}{m}{n}{<-> mathx10}{}
\DeclareSymbolFont{mathx}{U}{mathx}{m}{n}
\DeclareMathAccent{\widecheck}{0}{mathx}{"71}

\newenvironment{nouppercase}{%
  \renewcommand{\uppercasenonmath}[1]{}}{}

\begin{document}

\title{\Large Combinatorics of linear stability for Hamiltonian systems in arbitrary dimension \\ \vspace{0.5cm} \small On GIT quotients of the symplectic group, and the associahedron}

\author{Agustin Moreno}

\address[A.\ Moreno]{Universit\"at Heidelberg \\ Mathematisches Institut \\ Heidelberg \\ Germany}

\email{\href{mailto:agustin.moreno2191@gmail.com}{agustin.moreno2191@gmail.com}}

\author{Francesco Ruscelli}

\address[F.\ Ruscelli]{Universit\"at Heidelberg \\ Mathematisches Institut \\ Heidelberg \\ Germany}

\email{\href{mailto:fruscelli@mathi.uni-heidelberg.de}{fruscelli@mathi.uni-heidelberg.de}}




\begin{abstract}
    We address the general problem of studying linear stability and bifurcations of periodic orbits for Hamiltonian systems of arbitrary degrees of freedom. We study the topology of the GIT sequence introduced by the first author and Urs frauenfelder in \cite{FM}, in arbitrary dimension. In particular, we note that the combinatorics encoding the linear stability of periodic orbits is governed by a quotient of the associahedron. Our approach gives a topological/combinatorial proof of the classical Krein--Moser theorem, and refines it for the case of symmetric orbits.
\end{abstract}

\begin{nouppercase}
\maketitle    
\end{nouppercase}

\tableofcontents






\section{Introduction}

The stability of periodic orbits is a central topic in the study of Hamiltonian systems, going back to the problem of stability of the solar system in celestial mechanics. Ubiquitous in the study of ODEs, the notion of stability arises whenever studying orbits in families and their bifurcations, a practice which entails both theoretical and practical interest. For instance, from the perspective of space mission design, orbits used for parking a spacecraft around a target Moon should be as stable as possible, in order to minimize fuel corrections and station-keeping. From a mathematical point of view, the key notions of stability of a system come in three flavors, related by the following implications:
$$
\text{Non-linear (Lyapunov) stability }\Rightarrow \text{ linear stability } \Rightarrow \text{ spectral stability}.
$$
Non-linear stability, roughly speaking, means that trajectories which start near a given periodic orbit stay near the orbit for all time. Linear stability corresponds to stability of the origin for the linearized dynamics, i.e.\ the orbits of the linearized system should stay bounded. For a Hamiltonian system, this means that the eigenvalues of the monodromy matrix of the corresponding orbit should lie in the unit circle, and be semi-simple. Spectral stability, on the other hand, requires that eigenvalues all lie in the unit circle, but allows them to have multiplicity (so that orbits can escape to infinity in polynomial time, rather than exponential). In this paper, we will focus on the notion of \emph{linear} stability.

In the presence of symmetry, the study of linear stability of periodic orbits which are preserved by the symmetry can be significantly refined. With this end in mind, the first author and Urs Frauenfelder introduced in \cite{FM} the notion of the \emph{GIT sequence}, as a refinement of the \emph{Broucke stability diagram} \cite{Br69}, via the notion of \emph{B-signature}. The GIT sequence consists of a sequence of three spaces and maps between them whose topology encodes stability and bifurcations of periodic orbits, as well as their eigenvalue configurations, and provides obstructions to the existence of regular cylinder of orbits. In low dimensions, the spaces can be visualized in the plane or in three-dimensional space, which makes them amenable for numerical work. We should note that while the GIT sequence is designed to study linear stability, it blurs its distinction with spectral stability.

In \cite{FM}, the authors only dealt with the GIT sequence in the low dimensional cases. The aim of this paper is to completely determine the topology and combinatorics of the spaces in the sequence, in arbitrary dimension. In particular, we will determine their branching structure, and find normal form representatives, and therefore compute the degrees of the maps in the GIT sequence, which vary over the the components of a suitable decomposition of the base $\mathbb R^n$, each of them labelled by eigenvalue configurations. Among these components, there is a special one, the \emph{stable} component, which corresponds to stable periodic orbits. We will show that its combinatorics is governed by a quotient of the \emph{associahedron}. This will be a generalization to arbitrary dimension of the results in \cite{FM} (see also Howard--Mackay \cite{HM87}, Howard--Dullin \cite{HD98}). Our results can then be directly used to study the linear stability of periodic orbits for Hamiltonian systems in arbitrary degrees of freedom. We have also emphazised computatibility, in the sense that the structures we work with can be implemented numerically. Our topological approach also gives an alternative proof of the classical Krein--Moser theorem, reviewed in Appendix \ref{app:Krein}, and in fact gives a refinement for the case of symmetric orbits. 

Indeed, recall that the Krein--Moser theorem gives a criterion for when a \emph{Krein biurcation} may occur (i.e.\ two elliptic eigenvalues of the monodromy matrix come together and then bifurcate out of the circle). Our refinement gives a similar criterion for the situation when two \emph{hyperbolic} eigenvalues come together at a hyperbolic eigenvalue of multiplicity two and then become complex, but for the case of \emph{symmetric} orbits. We call such a transition a \emph{$\mathcal{HN}$-transition}, and the high-multiplicity eigenvalue, the \emph{transit} eigenvalue. Whether or not such a transition may occur is completely determined by the $B$-signature of the transit eigenvalue. Namely, the following result is a consequence of our topological study of the symplectic group.

\begin{theorem}\label{thm:HN} Consider a Hamiltonian with arbitrary degrees of freedom, admitting a symmetry. Let $t\mapsto \gamma_t, t\in [0,1],$ be a family of symmetric periodic orbits, undergoing an $\mathcal{HN}$-transition. Then the $B$-signature of the transit eigenvalue is indefinite.
\end{theorem}

The definition of $B$-signature will be given in Section \ref{sec:B_signature}, and the proof of this theorem is obtained in Appendix \ref{app:Krein}.

\medskip

\textbf{Acknowledgements.} The authors are grateful to Urs Frauenfelder, whose ideas inspired this paper. A.\ Moreno is currently supported by the Sonderforschungsbereich TRR 191 Symplectic Structures in Geometry, Algebra and Dynamics, funded by the DFG (Projektnummer 281071066 – TRR 191), and also by the DFG under Germany's Excellence Strategy EXC 2181/1 - 390900948 (the Heidelberg STRUCTURES Excellence Cluster).

\section{Preliminaries} In order to recall the definition of the GIT sequence, we need the following notion.

\begin{definition}[\textbf{GIT quotient}] Let $G$ be a group acting on a topological space $X$ by homeomorphisms. The \emph{GIT quotient} is the quotient space $X//G$ defined by the equivalence relation $x\sim y $ if the closures of the $G$-orbits of $x$ and $y$ intersect, endowed with the quotient topology.     
\end{definition}
The condition $x\sim y$ then means that there exists sequences $g_n,h_n \in G$ such that $\lim_n g_n \cdot x =\lim_n h_n \cdot y \in X$. The point is that the naive quotient space $X/G$ might not be Hausdorff in general, and $X//G$ always is, in a universal way. 

Consider now a symplectic manifold $(M,\omega)$, and an autonomous Hamiltonian system $H: M\rightarrow \mathbb R$. An \emph{involution} is a map $\rho:(M,\omega)\rightarrow (M,\omega)$ satisfying $\rho^2=id$, and it is anti-symplectic if $\rho^*\omega=-\omega$. Its \emph{fixed-point locus} is $\mathrm{Fix}(\rho)=\{x\in M:\rho(x)=x\}$, which is a Lagrangian submanifold of $M$. An anti-symplectic involution $\rho$ is a \emph{symmetry} of the system if $H\circ \rho=H.$ A periodic orbit $x$ is \emph{symmetric} with respect to an anti-symplectic involution $\rho$ if $\rho(x(-t))=x(t)$ for all $t$. The \emph{symmetric points} of the symmetric orbit $x$ are the two intersection points of $x$ with $\mathrm{Fix}(\rho)$, i.e.\
$$x\big(0\big),\,\,x\big(\tfrac{T}{2}\big) \in \mathrm{Fix}(\rho).$$
In particular, half of the symmetric periodic orbit is a Hamiltonian chord (i.e.\ trajectory) from
$\mathrm{Fix}(\rho)$ to itself. Hence we can think of a symmetric periodic orbit in two ways,
either as a closed string, or as an open string from the Lagrangian $\mathrm{Fix}(\rho)$ to itself.

The monodromy matrix of a symmetric orbit at a symmetric point is a \emph{Wonenburger} matrix, i.e.\ it satisfies
\begin{equation*}\label{symsymp}
M=M_{A,B,C}=\left(\begin{array}{cc}
A & B\\
C & A^t
\end{array}\right)\in Sp(2n),
\end{equation*}
where
\begin{equation}\label{eq:Wonenburger}
B=B^t,\quad C=C^t,\quad AB=BA^t,\quad
A^tC=CA,\quad A^2-BC=id,
\end{equation}
equations which ensure that $M$ is symplectic. The eigenvalues of $M$ are determined by those of the first block $A$ (see \cite{FM}):
\begin{itemize}
    \item If $\lambda$ is an eigenvalue of $M$ then its \emph{stability index} $a(\lambda)=\frac{1}{2}(\lambda + 1/\lambda)$ is an eigenvalue of $A$. 
    \item If $a$ is an eigenvalue of $A$ then $\lambda(a)=a+\sqrt{a^2-1}$ is an eigenvalue of $M$, for any choice of complex square root.
\end{itemize}

Note that in order to write the monodromy matrix in Wonenburger form, we implicitly chose a basis for $\mathrm{Fix}(\rho)$ at a symmetric point of the orbit (and extended it to a symplectic basis). A different choice of basis amounts to acting with an invertible matrix $R \in GL_n(\mathbb{R})$, via
\begin{equation*}\label{act}
R_*\big(A,B,C\big)=\Big(RAR^{-1},RBR^t,(R^t)^{-1}CR^{-1}\Big),
\end{equation*}
i.e.,\ $M_{A,B,C}$ is replaced by $M_{R_*(A,B,C)}$. We denote the space of Wonenburger matrices by
$$
Sp^\mathcal{I}(2n)=\{M_{A,B,C}: A,B,C \mbox{ satisfy } (\ref{eq:Wonenburger})\},
$$
which comes with the above action of $GL_n(\mathbb R)$. 

By a beautiful result of Wonenburger, every symplectic matrix $M\in Sp(2n)$ can be written as a product of two linear anti-symplectic involutions, i.e.\ $M=I_1 I_2$. From this, it is straightforward to derive the following fact (see \cite{FM}):

\begin{thm}[Wonenburger]\label{thm:Wonenburger}
    Every symplectic matrix $M\in Sp(2n)$ is symplectically conjugated to a Wonenburger matrix.
\end{thm}

In other words, the natural map
$$
Sp^\mathcal{I}(2n)//GL_n(\mathbb R) \rightarrow Sp(2n)//Sp(2n),
$$
$$
[M_{A,B,C}] \mapsto [M_{A,B,C}],
$$
is surjective.

In the presence of a symmetric periodic orbit, the above algebraic fact has a geometric interpretation: the monodromy matrix at each point of the orbit (a symplectic matrix) is symplectically conjugated via the linearized flow to the monodromy matrix at any of the symmetric points of the orbit (a Wonenburger matrix).

\section{The B-signature}\label{sec:B_signature} In this section, we discuss the notion of \emph{$B$-signature}, or \emph{$B$-signs}, introduced by the first author and Urs Frauenfelder in \cite{FM}, associated to a symmetric periodic orbit (see also \cite{Ay22}).

Consider a Wonenburger matrix
\begin{equation*}\label{symsymp}
M=M_{A,B,C}=\left(\begin{array}{cc}
A & B\\
C & A^t
\end{array}\right)\in Sp^\mathcal{I}(2n)\subset Sp(2n).
\end{equation*}

Let $\lambda,1/\lambda$ be an elliptic or hyperbolic pair of eigenvalues of $M$, which is simple (i.e.\ each eigenvalue is of multiplicity $1$). Then the stability index $a=a(\lambda)=\frac{1}{2}(\lambda+\frac{1}{\lambda})$ is an eigenvalue of $A$ (which has the same eigenvalues as $A^t)$, and is also simple. Let $v$ be an eigenvector of $A$ and $w$ be an eigenvector of $A^t$ with eigenvalue $a$. Assume that the pair $\lambda,1/\lambda$ is ordered so that $\lambda$ has positive imaginary part in the elliptic case, or has absolute value which is greater than $1$ in the hyperbolic case (i.e.\ $\lambda$ is the \emph{principal} eigenvalue).

\begin{definition}[\textbf{B/C-signs}] The \emph{B-sign} of $\lambda$ is the sign of $w^tBw$, and the \emph{C-sign} of $\lambda$ is the sign of $v^tCv$, denoted
$$
\mathrm{sign}_B(\lambda)=\mathrm{sign}(w^tBw)= \pm, 
$$
$$
\mathrm{sign}_C(\lambda)=\mathrm{sign}(v^tCv)= \pm. 
$$  
By definition, we let $\mathrm{sign}_{B/C}(1/\lambda)=-\mathrm{sign}_{B/C}(\lambda)$.
\end{definition}

We can think that we have attached a sign to the \emph{pair} $\lambda,1/\lambda$, ordered as above. Alternatively, we can attach these signs to the eigenvalues of $A$, as they correspond to the principal eigenvalues of $M$, after choosing the complex square root with positive imaginary part in the expression $\lambda=a+\sqrt{1-a^2}$. That the numbers $w^tBw$ and $v^tCv$ are non-vanishing follows from Equations (\ref{eq:Wonenburger}), and the fact that we are assuming the eigenvalues are elliptic or hyperbolic. It is easily checked that these definitions are independent of the eigenvectors $v,w$, and of the basis chosen. Therefore, to each symmetric point of a symmetric periodic orbit we have associated, for each simple elliptic or hyperbolic pairs of eigenvalues, a $B/C$-sign. For a given Wonenburger matrix, after ordering the simple and real eigenvalues of its $A$ block in strictly increasing order (which gives an order of the corresponding pairs of elliptic or hyperbolic pairs of eigenvalues of $M$), we obtain an ordered tuple of $B/C$ signs of the form $(\pm, \dots, \pm)$, respectively called the \emph{$B/C$-signature} of the matrix. Therefore a symmetric periodic orbit has \emph{two} such signatures, one for each of the two symmetric points. Moreover, the $C$-sign is completely determined by the $B$-sign, and viceversa, so they provide the same information. Namely, the $B$-sign agrees with the $C$-sign if $\lambda$ is hyperbolic, and they disagree if $\lambda$ is elliptic (this can be seen by inspecting the normal forms provided in \cite{FM}).

The following result illustrates the uses of these signs, as they give information on the type of periodic orbit.

\begin{thm}[\cite{FMb}]
A symmetric periodic orbit of a Hamiltonian system with two degrees of freedom is negative hyperbolic if and only if its two $B$-signs are different.
\end{thm}

The $B/C$-signs can be also defined in a straightforward way to the case where the eigenvalues are not necessarily simple. Indeed, let $\lambda$ be an elliptic or hyperbolic eigenvalue of $M\in Sp^{\mathcal{I}}(2n)$. Let $E_{\lambda}$ and $E^t_{\lambda}$ denote the $a(\lambda)$-eigenspaces of $A$ and $A^t$ respectively. We can then view $B$ as a bilinear form on $E^t_\lambda$, via
$$
B(v,w)=v^t B w,
$$
for $v,w \in E^t_\lambda$. Similarly, we can view $C$ as a bilinear form on $E_\lambda$.

\begin{definition}[\textbf{B/C-signature}] The \emph{B-signature} of $\lambda$ is the signature of $B\vert_{E^t_\lambda}$, and the \emph{C-signature} of $\lambda$ is the signature of $C\vert_{E_\lambda}$, denoted
$$
\mathrm{sign}_B(\lambda)=\mathrm{sign}(B\vert_{E^t_\lambda}), 
$$
$$
\mathrm{sign}_C(\lambda)=\mathrm{sign}(C\vert_{E_\lambda}). 
$$
We define the $B/C$-signature of $1/\lambda$ as the $B/C$-signature of $\lambda$.
\end{definition}

Recall that the signature of a non-degenerate bilinear form $G$ is the pair $(p,q)$, where $p$ is the dimension of a maximal subspace where $G$ is positive definite, and $q$ is the dimension of a maximal subspace where $G$ is negative definite. The fact that bilinear forms above are non-degenerate follows from Equations (\ref{eq:Wonenburger}) and ellipticity/hyperbolicity of the eigenvalues. Given a Wonenburger matrix, we order the real eigenvalues of $A$ in (non-strictly) increasing order, and this gives an ordered tuple $((p_1,q_1),\dots,(p_m,q_m))$ of $B/C$-signatures, which we call the $B/C-signature$ of the matrix. In the case where the eigenvalues are simple as above, one replaces $(1,0)$ with a $+$, and $(0,1)$ with a $-$ (as these are the only possibilities).

\section{GIT sequence: low dimensions}

We now discuss global topological methods in the study of periodic orbits, following the exposition in \cite{AFKM}. These methods encode: bifurcations; stability; eigenvalue configurations; obstructions to existence of regular families; and $B$-signs, in a visual and resource-efficient way. 

The main tool is the GIT sequence. This is a sequence of three branched spaces (or \emph{layers}), arranged into top, middle, and bottom (or base), together with two maps between them, which collapse certain branches together. Each branch is labelled by the $B$-signs. A symmetric orbit gives a point in the top layer, and an arbitrary orbit, in the middle layer. The base layer is $\mathbb R^n$ (the space of coefficients of the characteristic polynomial of the first block $A$ of $M_{A,B,C}$). Then a family of orbits gives a path in these spaces, so that their topology encodes valuable information, as it may sometimes enforce bifurcations, i.e.\ provide obstructions to the existence of regular families. The details are as follows. 

\medskip

The GIT sequence is the sequence of maps and spaces given by
$$
Sp^\mathcal{I}(2n)// GL_n(\mathbb R) \rightarrow Sp(2n)//Sp(2n) \rightarrow M_{n\times n}(\mathbb R)//GL_n(\mathbb R),
$$
$$
[M_{A,B,C}]\mapsto [M_{A,B,C}]\mapsto [A].
$$

We see that the above maps are well-defined, by checking that $M_{A,B,C}$ and $M_{R_*(A,B,C)}$ are symplectically conjugated. We are also implicitly using Theorem \ref{thm:Wonenburger} to define the second map, which is independent on the Wonenburger representative. Moreover, for the base of the sequence, we have the following nice fact (see \cite{FM} for a proof).

\begin{proposition}
    The base of the GIT sequence $M_{n\times n}(\mathbb R)//GL_n(\mathbb R)$ is homeomorphic to $\mathbb R^n$, where the homeomorphism maps a matrix to the coefficients of its characteristic polynomial, i.e.\
$$
M_{n\times n}(\mathbb R)//GL_n(\mathbb R)\rightarrow \mathbb R^n,
$$
    $$
    [A] \mapsto (c_{n-1},\dots,c_0),
    $$
    whenever
    $$
    p_A(t)=\det(A-t\cdot\mathds 1)=(-1)^nt^n+c_{n-1}t^{n-1}+\dots+c_0.
    $$
\end{proposition}

\begin{definition} In arbitrary dimension, we define the \emph{stability point} of the Wonenburger matrix $M_{A,B,C}\in Sp^{\mathcal I}(2n)$ as
$$
p=(s_1(\mu_1,\dots,\mu_n),\dots,s_n(\mu_1,\dots,\mu_n))\in \mathbb R^n,
$$
where $\mu_1,\dots,\mu_n$ are the eigenvalues of $A$, and $s_j$ is the $j$-th elementary symmetric polynomial, given by
$$
s_j(\mu_1,\dots,\mu_n)=\sum_{1\leq i_1 < \dots < i_j \leq n} \mu_{i_1}\dots\mu_{i_j}.
$$
\end{definition}

Then, by the above, the stability point is the result of applying the GIT sequence of maps to the given matrix.

\subsection{GIT sequence: 2D} We start with the simplest case, i.e.\ the case of an autonomous Hamiltonian of two degrees of freedom, so that the reduced monodromy matrix is an element in $Sp(2)=SL(2,\mathbb R)$.

Let $\lambda$ eigenvalue of $M\in Sp(2)$, with stability index $a(\lambda)=\frac{1}{2}(\lambda+1/\lambda)$. Then 
\begin{itemize}
    \item $\lambda=\pm 1$ iff $a(\lambda)=\pm 1$;
    \item $\lambda$ positive hyperbolic iff $a(\lambda)>1$; 
    \item $\lambda$ negative hyperbolic iff $a(\lambda)<-1$; and 
    \item $\lambda$ elliptic (i.e.\ stable) iff $-1<a(\lambda)<1$.
\end{itemize}
The Broucke stability diagram is then simply the real line, split into three components; see Figure \ref{fig:GIT_sequence}. If two orbits lie in different components of the diagram, then there are always bifurcations in any family joining them, as the topology of the diagram implies that any path between them has to cross the $\pm 1$ eigenvalues (corresponding respectively to bifurcation or period-doubling bifurcation).

One can think that the stability index ``collapses'' the two elliptic branches in the middle layer of Figure \ref{fig:GIT_sequence} together. These two branches are distinguished by the $B$-signs, coinciding with the Krein signs \cite{Kre2,Kre3}. There is an extra top layer for symmetric orbits, where now each hyperbolic branch separates into two, and there is a collapsing map from the top to middle layer. Note that to go from one branch to the other (say from the positive hyperbolic branch I to the positive hyperbolic branch II), the topology of the top layer implies that the eigenvalue 1 needs to be crossed. This means that one should expect bifurcations in any (symmetric) family joining them, \emph{even if} they project to the same component of the Broucke diagram. In this way, the information given by the diagram is much more refined for the case of symmetric orbits. If we say that two orbits are \emph{qualitatively equivalent} if they can be joined by a regular orbit cylinder, then the topology of the spaces in the GIT sequence give criteria to determine whenever two orbits are \emph{not} qualitatively equivalent. To sum up:

\begin{figure}[h]
    \centering
    \includegraphics[width=0.8\linewidth]{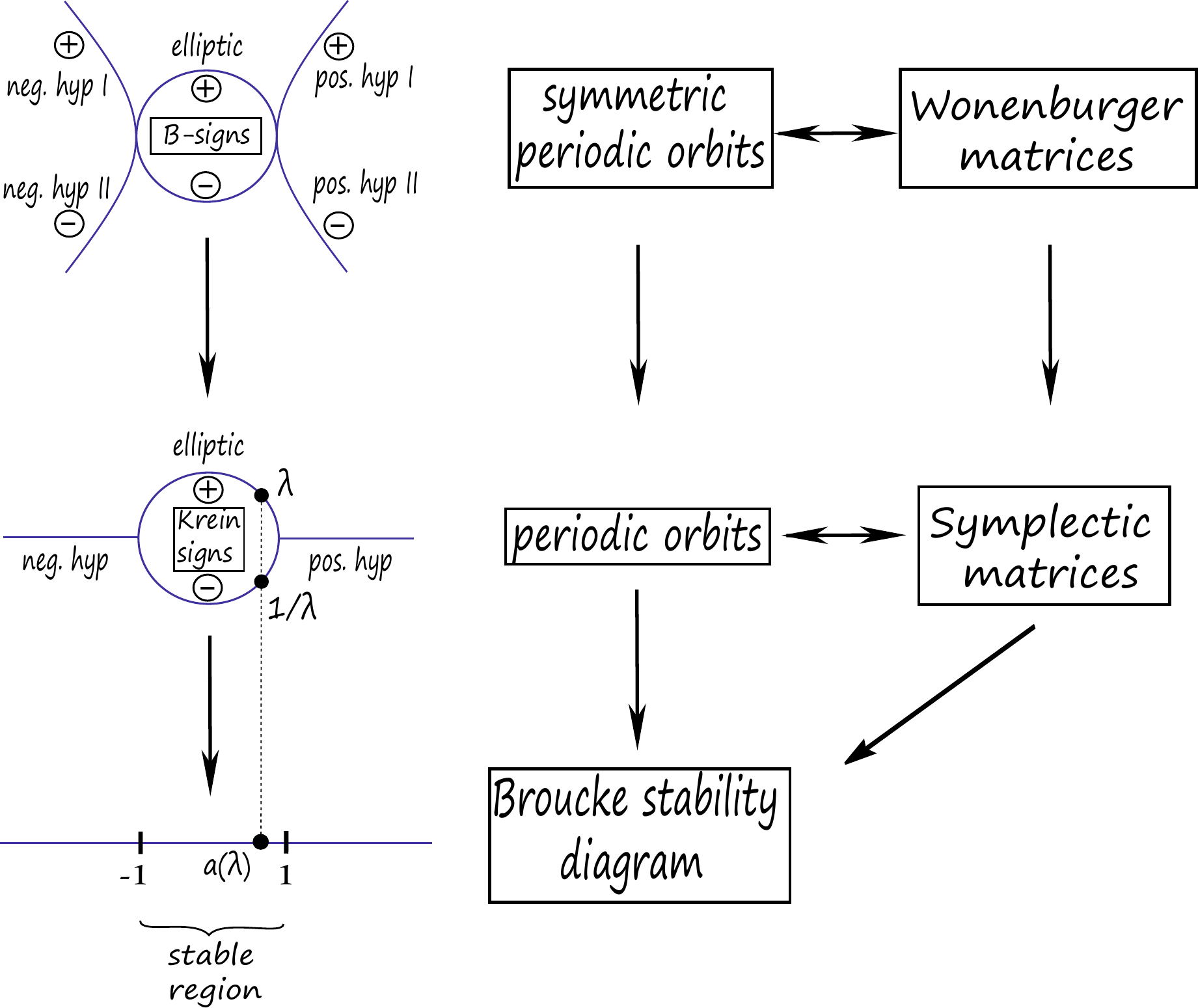}
    \caption{The 2D GIT sequence. One obtains more refined information for symmetric orbits.}
    \label{fig:GIT_sequence}
\end{figure}

\begin{itemize}
    \item $B$-signs ``separate'' hyperbolic branches, for symmetric orbits.

\medskip
    
    \item If two orbits lie in different components of the Broucke diagram, there are always bifurcations in any path joining them.

 \medskip
 
    \item If two symmetric orbits lie in the same component of the Broucke diagram, but if $B$-signs differ, one should \emph{also} expect bifurcation in any (symmetric) path joining them\footnote{We cautiosuly use the word ``\emph{expect}'' rather than give a mathematical statement, as theoretically orbits could tangentially pass through the Maslov cycle without birfurcating.}. 
    \end{itemize}

\subsection{GIT sequence: 3D} Now we apply the same idea, but for autonomous Hamiltonian systems with three degrees of freedom, for which reduced monodromy matrices are elements in $Sp(4)$.

\begin{figure}
    \centering
    \includegraphics[width=0.85\linewidth]{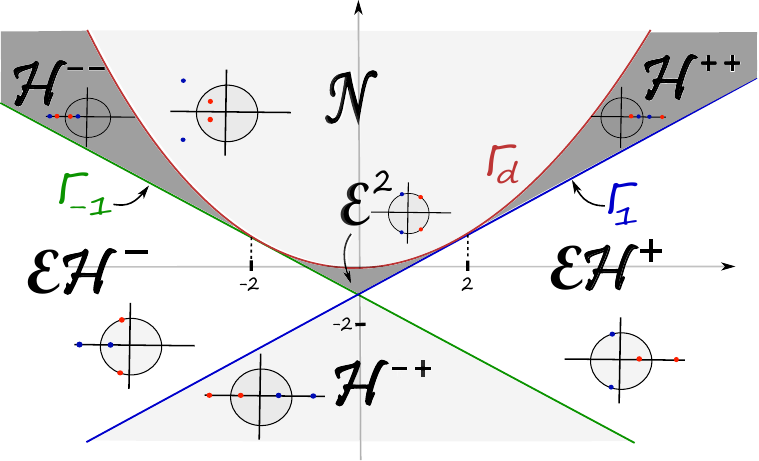}
    \caption{The 3D Broucke stability diagram. Here, $\Gamma_{\pm 1}$ corresponds to eigenvalue $\pm 1$, $\Gamma_d$ to double eigenvalue, $\mathcal{E}^2$ to doubly elliptic (the \emph{stable region}), and so on \cite{FM}.}
    \label{fig:Broucke_3D}
\end{figure}

Given a Wonenburger matrix $M=M_{A,B,C}\in Sp^\mathcal{I}(4)$, its \emph{stability point} is $$p=(\mbox{tr}(A),\det(A)) \in \mathbb R^2.$$ This point lies in the plane, which splits into regions corresponding to the eigenvalue configuration of $M$, as in Figure~\ref{fig:Broucke_3D}, which shows the Broucke stability diagram in $\mathbb R^2$. Each component is labelled according to the eigenvalue configuration. For instance, $\mathcal{E}^2$ (the \emph{doubly elliptic} component) corresponds to two pairs of elliptic eigenvalues; $\mathcal{EH}^{+}$ (the \emph{elliptic-positive hyperbolic} component), to one elliptic pair and a positive hyperbolic pair; $\mathcal{N}$, to complex quadruples, and so on. The parabola $\Gamma_d=\{y=x^2/4\}$ corresponds to double eigenvalues, i.e.\ two eigenvalues come together. The lines $\Gamma_{\pm 1}$ tangent to $\Gamma_d$ and with slope $\pm 1$ corresponds to matrices with eigenvalue $\pm 1$ in their spectrum. 

\begin{figure}
    \centering
    \includegraphics[width=1.13\linewidth]{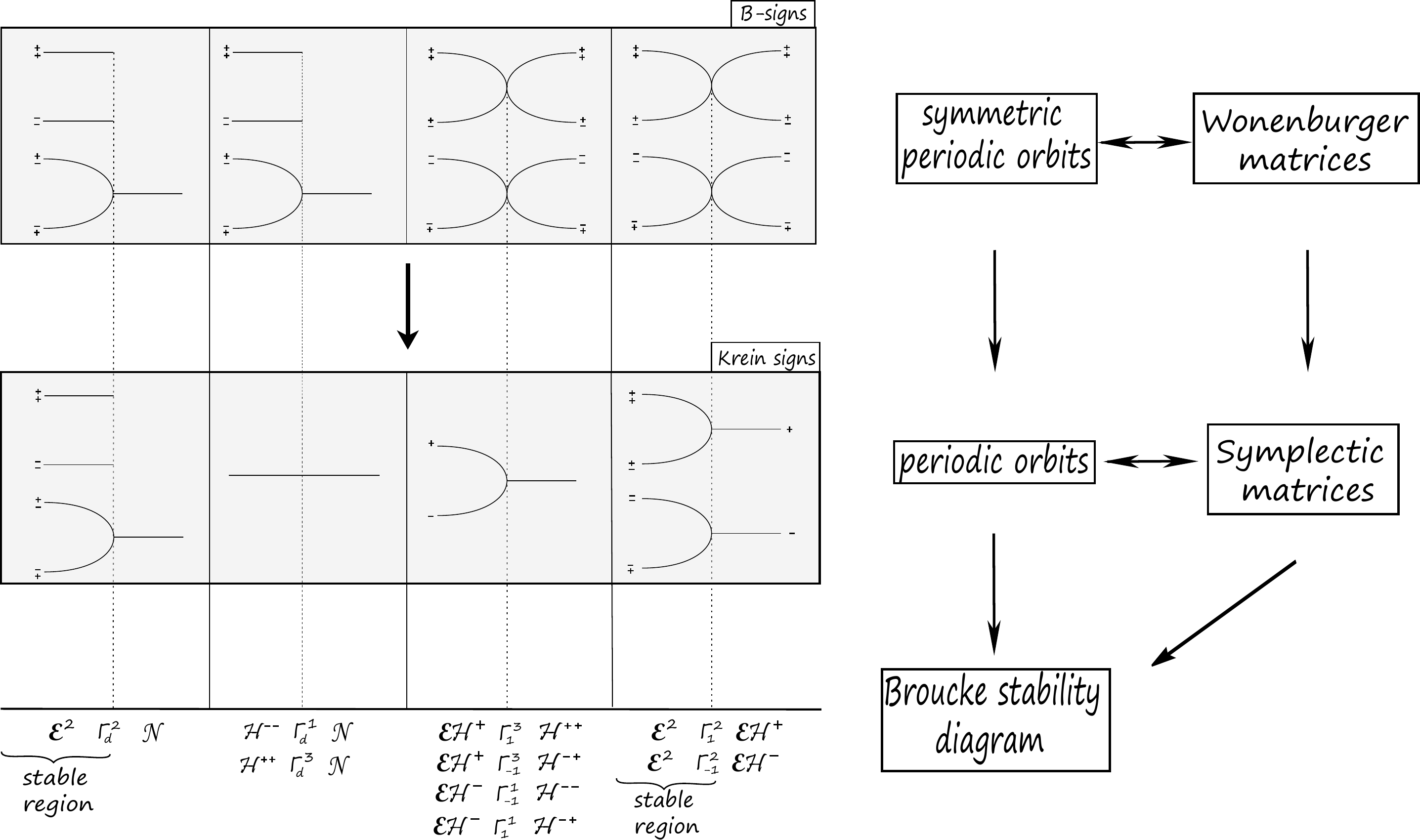}
    \caption{The branches (represented as lines) are two-dimensional, and come together at the 1-dimensional ``branching locus'' (represented as points), where we cross from one region to another of the Broucke diagram. The $1$-dimensional loci collapse to points over each of the three singular points $(2,1),(0,-1), (-2,1) \in \mathbb R^2$.}
    \label{fig:GIT_sequence_3D}
\end{figure}

The GIT sequence \cite{FM} adds two layers to this diagram, as shown in Figure~\ref{fig:GIT_sequence_3D}. The top layer has two extra branches than the middle one, for each hyperbolic eigenvalue. While the combinatorics and the global topology of the spaces involved is more complicated than the 2D case, the intuitive idea is still the same, i.e.\ that the amount of information for symmetric orbits is richer, and that we can distinguish more orbits up to qualitative equivalence. Note that as in this dimension we have two \emph{pairs} of eigenvalues, the $B$-signature is a pair $(\pm,\pm)$ of signs, and therefore the top layer has 4 branches over each component of the Broucke diagram (except the nonreal component).

\medskip

\begin{figure}
    \centering
\includegraphics[width=0.65\linewidth]{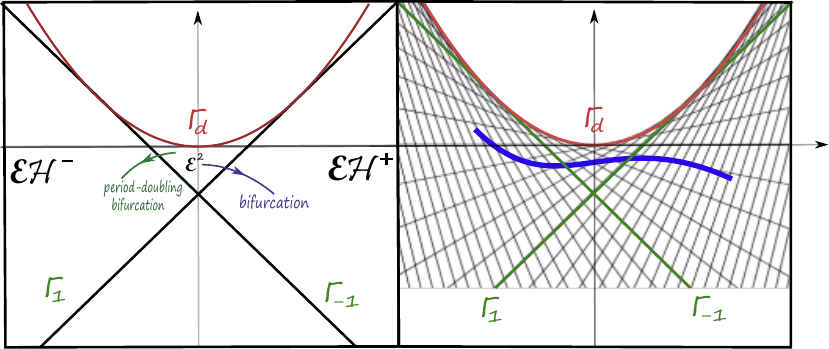}
\caption{Bifurcations are encoded by a pencil of lines.}
   \label{fig:complete_bifurcation}
\end{figure}

\textbf{Bifurcations in the Broucke diagram.} An orbit family $c\mapsto \gamma_c$ of symmetric orbits gives a path $c\mapsto p_c\in \mathbb R^2$ of stability points. The family bifurcates if $p_c$ crosses $\Gamma_1$. More generally, let $\Gamma^e_{\varphi}$ be the line with slope $\cos(2\pi \varphi)\in [-1,1]$ tangent to $\Gamma_d=\{y=x^2/4\}$, which corresponds to matrices with eigenvalue $e^{2\pi i\varphi}$; and $\Gamma^h_\lambda$ the tangent line with slope $a(\lambda)\in \mathbb R\backslash [-1,1]$, which corresponds to matrices with eigenvalue $\lambda$. Then a $k$-fold bifurcation happens when crossing $\Gamma^e_{l/k}$ for some $l$, i.e.\ the eigenvalue $e^{2\pi i l/k}$ is crossed. That is, higher order bifurcations are encoded by a pencil of lines tangent to a parabola, as in Figure~\ref{fig:complete_bifurcation}. Note that two such lines intersect at a point, which lies in a component determined by the lines (e.g.\ $\Gamma^e_\varphi \cap \Gamma^h_\lambda$ lies in $\mathcal{EH}^+$ if $\lambda>1$, and so on).

\section{GIT sequence: arbitrary dimension} We will now delve into the details of the GIT sequence, in arbitrary dimension. We will study the topology of the GIT quotients in the sequence, i.e.\ we will determine their branching structure. In particular, we will find normal form representatives in $Sp^{\mathcal{I}}(2n)//GL_n(\mathbb R)$, i.e.\ of Wonenburger matrices up to the natural action of $GL_n(\mathbb R)$, and therefore compute the degrees of the maps in the GIT sequence, which vary over the the components of a suitable decomposition of the base $\mathbb R^n$, each of them labelled by eigenvalue configurations. Among these components, there is a special one, the \emph{stable} component, which corresponds to stable periodic orbits. We will show that its combinatorics is governed by a quotient of the \emph{associahedron}. 

\medskip

\subsection{Some real algebraic geometry} Consider the space of monic polynomials with real coefficients and of degree $n$, i.e.\ of the form 
$$
p(t)=(-1)^nt^n +c_{n-1}t^{n-1}+\dots + c_0,
$$
with $c_i\in \mathbb R$, which we identify with $\mathbb R^n$ via
$$
p \longleftrightarrow (c_0,\dots,c_{n-1}).
$$
Recall that the \emph{discriminant} of a polynomial is defined as the expression
$$
\Delta(p)=\prod_{i<j} (\lambda_i-\lambda_j)^2,
$$
where $\lambda_1\dots,\lambda_n$ are the complex roots of $p$. Then $p$ has a multiple root if and only if $\Delta(p)=0$, and if these roots are all real and simple, then $\Delta(p)>0$. Moreover, as $\Delta(p)$ is a symmetric polynomial in the roots of $p$, it follows from the fundamental theorem of symmetric polynomials that it is a polynomial in the basic symmetric polynomials of the roots, i.e.\ a polynomial in the coefficients $c_i$ of $p$. Therefore 
$$
V(\Delta)=\{\Delta=0\}\subset \mathbb R^n
$$
is a (in general, singular) real algebraic variety, the \emph{discriminant variety}. For instance, the discriminant of a degree $2$ polynomial $p(x)=x^2+bx+c$ is $\Delta(p)=b^2-4c$, so that $V(\Delta)=\{b^2=4c\}$ is a parabola in $\mathbb R^2$, and that of a degree $2$ polynomial $p(x)= -x^3+bx^2+cx+d$ is $\Delta(p)=b^2c^2+4c^3-4b^3d-27d^2-18bcd,$ so that $V(\Delta)$ is a quadric hypersurface in $\mathbb R^3$. 

For $a\in \mathbb R$, we let $\Gamma_a=\{p: p(a)=0\}\subset \mathbb R^n$. Note that the equation $p(a)=0$ is linear in the coefficients of $p$, and so $\Gamma_a$ is a linear hyperplane in $\mathbb R^n$, given by $\Gamma_a=\{P_a=0\}$ for a linear equation $P_a$. Similarly, for $\alpha\in \mathbb C\backslash \mathbb R$, we let $V_\alpha=\{p: p(\alpha)=0\}\subset \mathbb R^n$. Note that $V_\alpha=V_{\overline{\alpha}}$, as complex roots come in conjugate pairs. We then define the regions of $\mathbb R^n$ given by
$$
\mathcal{E}=\bigcup_{a\in[-1,1]}\Gamma_a,\; \mathcal{H}^+=\bigcup_{a\geq 1}\Gamma_a,\; \mathcal{H}^-=\bigcup_{a\leq -1}\Gamma_a, \;\mathcal{N}=\bigcup_{\alpha \in \mathbb C\backslash \mathbb R} V_\alpha,
$$
respectively the \emph{elliptic, positive/negative hyperbolic, nonreal} regions. Note that the intersection $\Gamma_{aa'}=\Gamma_a\cap \Gamma_{a'}$ of two hyperplanes $\Gamma_a,\Gamma_{a'}$ for $a\neq a'$ is a codimension-$2$ affine linear subspace of $\mathbb R^n$ corresponding to polynomials $$p(x)=(x-a)(x-a')q(x)$$ for some $q\in \mathbb R^{n-2}$. Similarly, $\Gamma_{a_1\dots a_n}=\Gamma_{a_1}\cap \dots \cap \Gamma_{a_n}$ is a point, corresponding to the polynomial
$$
p(x)=(-1)^n(x-a_1)\dots (x-a_n).
$$
Similarly, $V_{\alpha_1\dots \alpha_n}=V_{\alpha_1}\cap\dots\cap V_{\alpha_n}$ is a point. If we let $a'\rightarrow a$, we see that $\Gamma_{aa'}$ converges to $\Gamma_a$, but also becomes tangent to the discriminant variety along the subspace $$V_a(\Delta)=\{(x-a)^2q(x): q \in \mathbb R^{n-2}\}=V(\Delta)\cap \Gamma_a\subset V(\Delta).$$ 
Note that the \emph{regular} points of $V(\Delta)$ are those points $(x-a)^2q(x)$ which satisfy $q(a)\neq 0$, and so $V(\Delta)$ has a well-defined tangent space at these points. By comparing dimensions, we see that $\Gamma_a$ coincides with this tangent space at each of the regular points. Moreover, $\Gamma_a$ is tangent to $V(\Delta)$ along $V_a(\Delta)$. Varying $a$, we further see that the discriminant variety is the \emph{envelope} of the family of hyperplanes $\{\Gamma_a\}_{a\in \mathbb R}$, i.e.\ it is the locus of points $V(\Delta)=\bigcup_a V_a(\Delta)$ which is tangent to each $\Gamma_a$, and is the union of these tangency loci (the \emph{characteristic subspaces} $V_a(\Delta)$). The envelope $V(\Delta)$ can then be described as the set of coefficients $c\in \mathbb R^n$ satisfying
$$
P_a(c)=\partial_aP_a(c)=0
$$
for some $a$.

We also have a decomposition of $\mathbb R^n$ into regions labelled by the root configurations of the corresponding polynomials. Namely, for $k,l,m,r$ with $k+l+m+2r=n$, we let
$$
\mathcal{M}^{klmr}_0:=(\mathcal{H}^-)^k\mathcal{E}^l(\mathcal{H}^+)^m\mathcal{N}^r:=\bigcup_{\substack{a_1,\dots,a_l\in [-1,1]\\ b_1,\dots, b_m\geq 1 \\ c_1,\dots, c_k\leq -1\\ \alpha_1,\dots, \alpha_r\in \mathbb C\backslash\mathbb R}}\Gamma_{a_1\dots a_l}\cap \Gamma_{b_1\dots b_m} \cap \Gamma_{c_1 \dots c_k} \cap V_{\alpha_1\dots\alpha_r}.
$$
Then 
$$
\mathbb R^n=\bigcup_{\substack{k,l,m,r\geq 0\\ k+l+m+2r=n}}(\mathcal{H}^-)^k\mathcal{E}^l(\mathcal{H}^+)^m\mathcal{N}^r.
$$
See Figure \ref{fig:Broucke_3D} for the case $n=2$ (where $V(\Delta)$ is denoted $\Gamma_d$). 

Among the regions in this decomposition, there is a unique compact one, corresponding to
$$
\mathcal{E}^n:=(\mathcal{H}^-)^0\mathcal{E}^n(\mathcal{H}^+)^0\mathcal{N}^0.
$$
This is the \emph{stable} region. We will give a combinatorial description of this component below.

\begin{figure}
    \centering
    \includegraphics[width=0.8\linewidth]{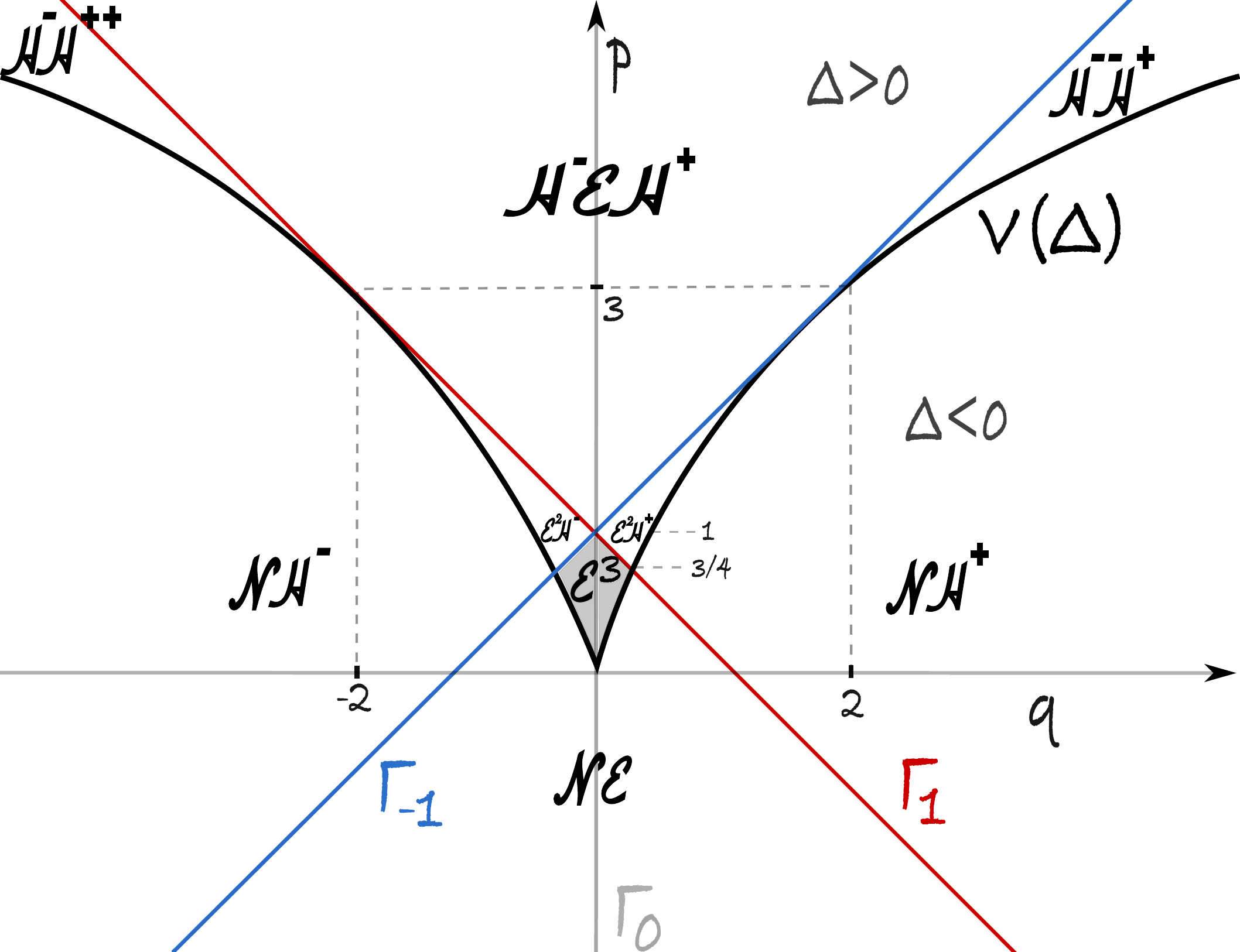}
    \caption{The stability diagram for depressed cubics.}
    \label{fig:depressed_diagram}
\end{figure}

\begin{example} For the case $n=3$, every polynomial
$$
A(x)=-x^3+bx^2+cx+d
$$  
can be transformed via the change of variables $y=x-b/3$ to a polynomial with no degree $2$ term (a \emph{depressed} cubic polynomial), i.e.\ of the form
$$
B(y)=-y^3+py+q
$$
with $p=-c-\frac{b^2}{3},\; q=-\frac{2b^3+9bc}{27}-d$. Its discriminant is now
$$
\Delta(B)=4p^3-27q^2.
$$
If we identify the space of depressed cubics with $\mathbb R^2$ via $B\leftrightarrow (q,p)$, the vanishing locus of $\Delta$ is a cubic in the plane singular at the origin (corresponding to the polynomial $y^3$, having triple root), and the decomposition of the plane according to the root configuration is depicted in Figure \ref{fig:depressed_diagram}. Note that there are no $(\mathcal{H}^{+})^3$ or $(\mathcal{H}^{-})^3$ components, as no depressed cubic can have roots which are all greater than $1$ in absolute value (as they sum to zero). The pencil of lines $\{\Gamma_a\}_{a\in \mathbb R}$ has $V(\Delta)$ as envelope, with $\Gamma_a=\{P_a=-a^3+pa+q=0\}$ of slope $-1/a$, corresponding to polynomials with root $a$. The discriminant variety is the loci of points satisfying 
$$
P_a=-a^3+pa+q=\partial_aP_a=-3a^2+p=0
$$
for some $a$, as is easily checked. Note moreover that if two lines $\Gamma_a$ and $\Gamma_{a'}$ are given, then there is a unique third line $\Gamma_{a''}$ that passes through the intersection $\Gamma_{aa'}$. In other words, if we impose that we get a depressed cubic, then the third root $a''$ is determined by the two other roots $a,a'$ via the equation $a+a'+a''=0$. 

\end{example}

\textbf{Regular cases.} We now study the GIT sequence. We first consider the regular cases, i.e.\ when the matrix has only simple eigenvalues different from $\pm 1$. Denote by $\mathcal{M}^{klmr}$ the stratum of matrices $[M]\in Sp^\mathcal{I}(2n)//GL_n(\mathbb R)$ with $M=M_{A,B,C}$, the $A$-block having an eigenvalue configuration with
$$
\mu^-_1\leq \dots \leq \mu^-_k\leq -1\leq \nu_{1}\leq \dots \leq \nu_l \leq 1 \leq \mu^+_{1}\leq\dots\leq \mu_m^+ 
$$
the real eigenvalues of $A$, and 
$$
\alpha_1,\overline{\alpha}_1,\dots, \alpha_r,\overline{\alpha}_r
$$
the complex conjugate eigenvalues of $A$, satisfying
$$
k+l+m+2r=n.
$$

Then there exist unique $\beta^\pm_j, r_j \in [0,+\infty)$ and $\theta_j,\gamma_j\in [0,\pi]$ such that
$$
\mu_j^-=-\cosh \beta_j^-,\;\mu_j^+=\cosh \beta_j^+,\; \nu_j=\cos \theta_j,\;\alpha_j=r_je^{i\gamma_j},\;\overline{\alpha}_j=r_je^{-i\gamma_j}. 
$$ 
Moreover, up to conjugation $A$ can be assumed to split as a sum
$$
A=\mathrm{diag}\left(\mu^-_1,\dots,\mu^-_k,\nu_{1},\dots,\nu_l,\mu^+_{1},\dots,\mu_m^+\right)\oplus\bigoplus_{j=1}^rR(r_j,\gamma_j), 
$$
where 
$$
R(r,\gamma)=\left(\begin{array}{cc}
   r\cos \gamma & -r\sin \gamma \\
   r\sin \gamma  & r\cos \gamma
\end{array} \right)
$$
is the composition of a dilation with a rotation. Therefore each summand can be treated separately, which is what is done in \cite{FM}. The result is the normal form
\begin{equation}
    \begin{split}
B= &\; \mathrm{ diag}\left(\epsilon^-_1\sinh(\beta_1^-),\dots,\epsilon^-_k\sinh(\beta_k^-),\delta_1\sin \theta_j,\dots,\delta_l\sin \theta_l,\epsilon^+_1\sinh(\beta_1^+),\dots,\epsilon^+_m\sinh(\beta_m^+)\right) \\
&\oplus\bigoplus_{j=1}^r\mathrm{diag}(1,-1),\\
C= &\; \mathrm{ diag}\left(\epsilon^-_1\sinh(\beta_1^-),\dots,\epsilon^-_k\sinh(\beta_k^-),-\delta_1\sin \theta_j,\dots,-\delta_l\sin \theta_l,\epsilon^+_1\sinh(\beta_1^+),\dots,\epsilon^+_m\sinh(\beta_m^+)\right)\\ 
&\oplus\bigoplus_{j=1}^rS(r_j,\gamma_j),\\ 
    \end{split}
\end{equation}
where $\epsilon_j^\pm,\delta_j$ are $\pm$ signs, and
$$
S(r,\gamma)=\left(\begin{array}{cc}
r^2 \cos 2\gamma -1 & -r^2\sin 2\gamma\\
-r^2 \sin 2\gamma & -r^2 \cos 2\gamma +1
\end{array}\right).
$$
We then see that the $B$-signature of $M$ is
$$
\mathrm{sign}_B(M)=\left(\epsilon^-_1,\dots,\epsilon^-_k,\delta_1,\dots,\delta_l,\epsilon^+_1,\dots,\epsilon^+_m\right),
$$
so that they are distinct as elements in the GIT quotient $Sp^\mathcal{I}(2n)//GL_n(\mathbb R)$. Note that if we flip some of the $\epsilon^{\pm}_i$ signs, we obtain matrices that are symplectically conjugated to the original one. 

\begin{remark}\label{rk:Krein-Moser}
    Note that if there are complex quadruples in the spectrum, then the $B$-block always has at least one summand of the form diag$(1,-1)$. Viewed as a bilinear form, this matrix always has \emph{mixed} signature. As signatures behave continuously in the space of non-degenerate bilinear forms, this implies that the corresponding quadruple cannot be connected to a hyperbolic or elliptic pair of multiplicity two of \emph{definite} signature, while fixing the remaining eigenvalues. This is the main observation that implies the Krein--Moser theorem, cf.\ Appendix \ref{app:Krein}. This is also what implies its refinement for symmetric orbits (Theorem \ref{thm:HN}).
\end{remark}

By varying the $B$-signature we see that there are $2^{n-2r}=2^{k+l+m}$ connected components of the interior of the loci $\mathcal{M}^{klmr}$, i.e.\
$$
\mathrm{int}\;\mathcal{M}^{klmr}=\bigsqcup_{\epsilon}\mathcal{M}^{klmr; \epsilon},
$$
where $\epsilon$ is a $B$-signature with $n-2r$ entries. Moreover, if we denote by $\mathcal{M}^{klmr}_0$ the region of $\mathbb R^n$ which is the image of $\mathcal{M}^{klmr}$ under the GIT map $Sp^\mathcal{I}(2n)//GL_n(\mathbb R)\rightarrow \mathbb R^n$, then the covering degree of this map over $\mathcal{M}^{klmr}_0$ is precisely the number of connected components (or \emph{branches}) of $\mathcal{M}^{klmr}$, i.e.\
$$
\mathrm{deg}\;\mathcal{M}^{klmr}=2^{n-2r}=2^{k+l+m}.
$$
The GIT map simply collapses the different $\mathcal{M}^{klmr; \epsilon}$ to the region $\mathcal{M}^{klmr}_0$.

We can also consider the stratum $\mathcal{M}_1^{klmr}$ of matrices $[M]\in Sp(2n)//Sp(2n)$ with the above eigenvalue configuration. Then by varying the Krein-signature $\kappa$, we obtain a splitting
$$
\mathrm{int}\;\mathcal{M}_1^{klmr}=\bigsqcup_{\kappa}\mathcal{M}_1^{klmr; \kappa},
$$
into connected components, each of which is $n$-dimensional. The GIT sequence map $Sp^\mathcal{I}(2n)//GL_n(\mathbb R)\rightarrow Sp(2n)//Sp(2n)$ collapses certain branches $\mathcal{M}_1^{klmr; \kappa}$ together, whenever two such signatures $\kappa,\kappa'$ coincide up to removing all those corresponding to hyperbolic eigenvalues. More precisely, we define a \emph{collapsing} operation on signatures, given by
$$
\mathrm{Col}\left(\epsilon^-_1,\dots,\epsilon^-_k,\delta_1,\dots,\delta_l,\epsilon^+_1,\dots,\epsilon^+_k\right)=\left(\delta_1,\dots,\delta_l\right),
$$
which only keeps the signs of the elliptic eigenvalues. Then the GIT sequence sends
$$
\mathcal{M}^{klmr;\epsilon} \mapsto \mathcal{M}_1^{klmr;\mathrm{Col}(\epsilon)} \mapsto \mathcal{M}_0^{klmr}=(\mathcal{H}^-)^k\mathcal{E}^l(\mathcal{H}^+)^m\mathcal{N}^r
$$
homeomorphically. The degree of $\mathcal{M}_1^{klmr}$, defined as the covering degree of the second map in the GIT sequence over $\mathcal{M}_0^{klmr}$ (i.e.\ the number of branches $\mathcal{M}_1^{klmr;\kappa}$ that project to it), is then
$$
\mathrm{deg}\;\mathcal{M}_1^{klmr}=2^{l}.
$$
The degree of the first map in the GIT sequence over a branch $\mathcal{M}_1^{klmr;\kappa}$, the \emph{relative} degree (i.e.\ the number of branches $\mathcal{M}^{klmr;\epsilon}$ that collapse to it), is then
$$
\mathrm{deg}\;\left(\mathcal{M}^{klmr}\Big\vert\mathcal{M}_1^{klmr}\right)=2^{k+m},
$$
which is independent of $\kappa$. We have the obvious multiplicative formula
$$ \mathrm{deg}\;\left(\mathcal{M}^{klmr}\Big\vert\mathcal{M}_1^{klmr}\right)\cdot \mathrm{deg}\;\mathcal{M}_1^{klmr}=\mathrm{deg}\;\mathcal{M}^{klmr}.
$$
\medskip

\textbf{Nonregular cases.} The nonregular cases can be dealt with similarly, although the combinatorics gets more involved. Indeed, assume that $A$ has real eigenvalues 
$$
\mu^-_1\leq\dots \leq \mu^-_k\leq -1\leq \nu_{1}\leq \dots \leq \nu_l \leq 1 \leq \mu^+_{1}\leq \dots\leq \mu_m^+, 
$$
where we also allow $\pm 1$ as an eigenvalue, and complex eigenvalues
$$
\alpha_1,\overline{\alpha}_1,\dots,\alpha_r,\overline{\alpha}_r.
$$
We denote the multiplicities by
$$
m_j^\pm=\mathrm{mult}(\mu_j^\pm),\; o_j=\mathrm{mult}(\nu_j), \;M_\pm= \mathrm{mult}(\pm 1),\; p_j= \mathrm{mult}(\alpha_j)=\mathrm{mult}(\overline{\alpha_j}).
$$
If we let $m^-=(m_1^-,\dots,m_k^-), m^+=(m_1^+,\dots,m_m^+), o=(o_1,\dots,o_l), p=(p_1,\dots,p_r)$, we have
$$
\vert m^+ \vert + M_-+\vert o \vert + M_+ +\vert m^- \vert +2\vert p\vert=n, 
$$
where 
$$
\vert(a_1,\dots,a_n)\vert=\sum_{j=1}^na_i.
$$
We further denote by 
$$\mathcal{M}^{klmr}_{m^-M_-oM_+m^+p}=(\mathcal{H}^-)_{m^-}^k(-\mathcal I)_{M_-}\mathcal{E}_o^l\;\mathcal{I}_{M_+}(\mathcal{H}^+)_{m^+}^m\mathcal{N}_p^r$$ the stratum of matrices $[M]\in Sp^{\mathcal{I}}(2n)//GL_n(\mathbb R)$ with $A$-block having eigenvalue configuration as above.

\begin{figure}
    \centering
    \includegraphics[width=0.58\linewidth]{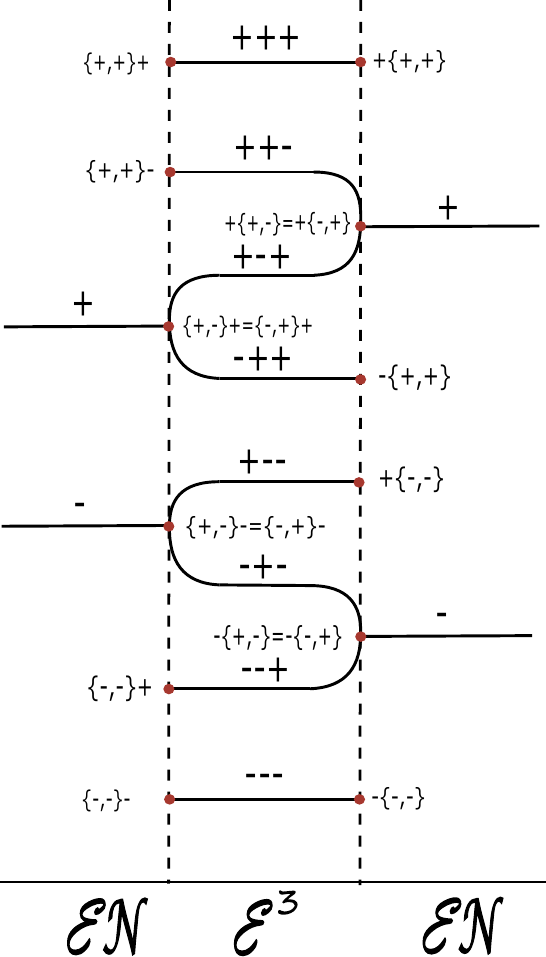}
    \caption{The branching structure of $Sp^{\mathcal{ I}}(6)//GL_3(\mathbb R)$, for the case $n=3$, corresponding to a transition from $\mathcal{E}^3$ to $\mathcal{EN}$. The right and left hand sides should be identified. The branching structure for $Sp(6)//Sp(6)$ coincides for this transition, as no hyperbolic eigenvalues are involved. Note that the branches do \emph{not} cross along a contraction which has definite signature inside the brackets, as predicted by the Krein--Moser theorem; see Appendix \ref{app:Krein}. This holds in all dimension (see Remark \ref{rk:Krein-Moser}), and gives another topological proof of this theorem.}
    \label{fig:EN}
\end{figure}

\begin{figure}
    \centering
    \includegraphics[width=0.8\linewidth]{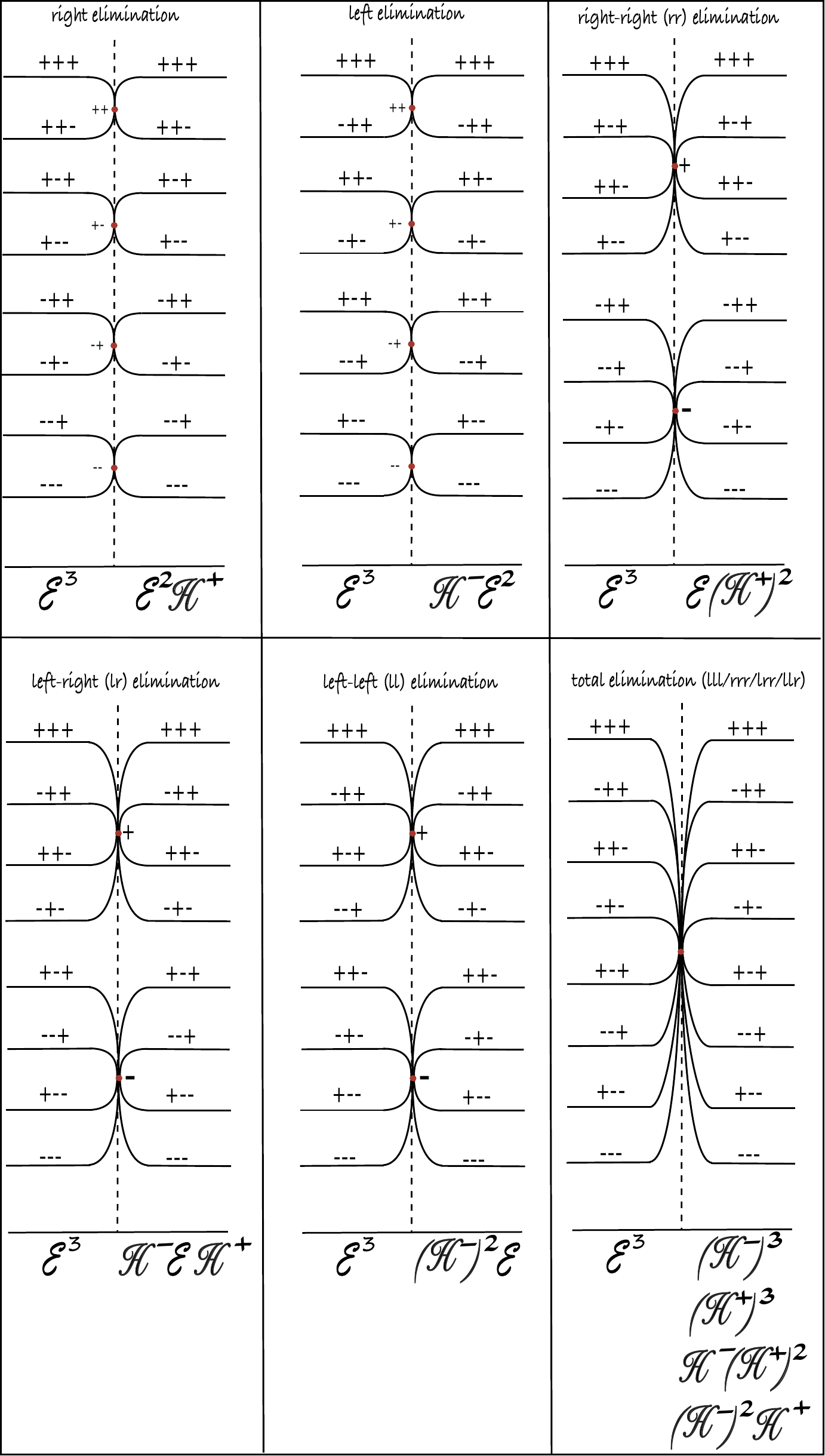}
    \caption{Case $n=3$, for $Sp^{\mathcal{I}}(6)//GL_3(\mathbb R)$. The transitions out of the stable region which do not involve $\mathcal{N}$ correspond to all possible elimination contractions.}
    \label{fig:eliminations}
\end{figure}

Moreover, there exist unique $\beta^\pm_j, r_j \in [0,+\infty)$ and $\theta_j,\gamma_j\in [0,\pi]$ such that
$$
\mu_j^-=-\cosh \beta_j^-,\;\mu_j^+=\cosh \beta_j^+,\; \nu_j=\cos \theta_j,\;\alpha_j=r_je^{i\gamma_j},\;\overline{\alpha}_j=r_je^{-i\gamma_j}. 
$$ 
As we can ignore Jordan blocks in the GIT quotient, we can assume that 
$$
A=\bigoplus_{j=1}^k \mu_j^- \mathds 1_{m_j^-} \oplus -\mathds 1_{M_-} \oplus \bigoplus_{j=1}^l \nu_j \mathds 1_{o_j} \oplus \mathds 1_{M_+} \oplus \bigoplus_{j=1}^m \mu_j^+ \mathds 1_{m_j^+} \oplus \bigoplus_{j=1}^r \bigoplus_{i=1}^{p_j} R(r_j,\gamma_j), 
$$
where $\mathds 1_n$ is the identity matrix of size $n$. Applying \cite{FM} to each summand, we obtain the normal form
\begin{equation}
\begin{split}
B=&\bigoplus_{j=1}^k \left(\sinh(\beta_j^-) \mathds 1_{a_j^-}  \oplus -\sinh(\beta_j^-) \mathds 1_{b_j^-}\right)\oplus 0_{M_-} \oplus \bigoplus_{j=1}^l \left(\sin(\theta_j) \mathds 1_{c_j}  \oplus -\sin(\theta_j) \mathds 1_{d_j}\right) \oplus 0_{M_+}\oplus\\
& \bigoplus_{j=1}^m  \left(\sinh(\beta_j^+) \mathds 1_{a_j^+}  \oplus -\sinh(\beta_j^+)\mathds 1_{b_j^+}\right) \oplus \bigoplus_{j=1}^r \bigoplus_{i=1}^{p_j} \mathrm{diag}(1,-1), 
\end{split}
\end{equation}
\begin{equation}
\begin{split}
C=&\bigoplus_{j=1}^k \left(\sinh(\beta_j^-) \mathds 1_{a_j^-}  \oplus -\sinh(\beta_j^-) \mathds 1_{b_j^-}\right)\oplus 0_{M_-} \oplus \bigoplus_{j=1}^l \left(-\sin(\theta_j) \mathds 1_{c_j}  \oplus \sin(\theta_j) \mathds 1_{d_j}\right) \oplus 0_{M_+}\oplus\\
& \bigoplus_{j=1}^m  \left(\sinh(\beta_j^+) \mathds 1_{a_j^+}  \oplus -\sinh(\beta_j^+)\mathds 1_{b_j^+}\right) \oplus \bigoplus_{j=1}^r \bigoplus_{i=1}^{p_j} S(r_j,\gamma_j), 
\end{split}
\end{equation}
with $0_m$ the zero matrix of size $m$, and with
$$
a_j^\pm+b_j^\pm=m_j^\pm,\; c_j+d_j=o_j.
$$
The $B$-signature of $M$ is then
$$
\mathrm{sign}_B(M)=((a_1^-,b_1^-),\dots,(a_k^-,b_k^-), (c_1,d_1),\dots,(c_l,d_l),(a_1^+,b_1^+),\dots,(a_m^+,b_m^+)).
$$
If we denote $\mathcal{M}^{klmr}_{m^-M_-oM_+m^+p;0}\subset \mathbb R^n$ the image of $\mathcal{M}^{klmr}_{m^-M_-oM_+m^+p}$ under the maps in the GIT sequence, by counting the normal forms for varying $B$-signature, we see that the degree of the GIT map over this component (i.e.\ the number of connected components $\mathcal{M}^{klmr;\epsilon}_{m^-M_-oM_+m^+p}$ or the \emph{branches} of $\mathcal{M}^{klmr}_{m^-M_-oM_+m^+p}$) is
$$
\mathrm{deg}\;\mathcal{M}^{klmr}_{m^-M_-oM_+m^+p}=\prod_{j=1}^k(m_j^-+1)\prod_{j=1}^l (o_j+1)\prod_{j=1}^m(m_j^++1).
$$
Note that if $m_j^\pm=o_j=1$ then we indeed recover the formula for the regular case. And if $m_j^\pm=o_j=0$ (i.e.\ only complex quadruples or $\pm 1$ as eigenvalues) then the above degree is $1$. The reader is invited to compare the above formula to Figures \ref{fig:GIT_sequence_3D}, \ref{fig:EN}, \ref{fig:eliminations}, as it gives the number of branches above each component of the base. Also, note that the dimension of the interior of each branch is
\begin{equation}
\begin{split}
\dim \text{int}\;\mathcal{M}^{klmr;\epsilon}_{m^-M_-oM_+m^+p}=&n-M_--M_+-\sum_{j=1}^k(m^-_j-1)-\sum_{j=1}^l(o_j-1)\\&-\sum_{j=1}^m(m^+_j-1)-2\sum_{j=1}^r(p_j-1)\\
&=k+l+m+2r\leq n.
\end{split}
\end{equation}

Similarly, we denote by $\mathcal{M}^{klmr}_{m^-M_-oM_+m^+p;1}$ the stratum of matrices $[M]\in Sp(2n)//Sp(2n)$ having eigenvalue configuration as above, and by $\mathcal{M}^{klmr;\kappa}_{m^-M_-oM_+m^+p;1}$ its connected components, labelled by the Krein signature $\kappa$. By counting these connected components, we see that its degree is
$$
\mathrm{deg}\;\mathcal{M}^{klmr}_{m^-M_-oM_+m^+p;1}=\prod_{j=1}^l (o_j+1).
$$
We also have a collapsing map
$$
\mathrm{Col}((a_1^-,b_1^-),\dots,(a_k^-,b_k^-), (c_1,d_1),\dots,(c_l,d_l),(a_1^+,b_1^+),\dots,(a_m^+,b_m^+))=((c_1,d_1),\dots,(c_l,d_l)),
$$
and the GIT sequence maps
$$
\mathcal{M}^{klmr;\epsilon}_{m^-M_-oM_+m^+p}\mapsto\mathcal{M}^{klmr; \mathrm{Col}(\epsilon)}_{m^-M_-oM_+m^+p;1}\mapsto \mathcal{M}^{klmr}_{m^-M_-oM_+m^+p;0},
$$
homeomorphically, and so their interiors have the same dimension. The relative degree is then
$$
\mathrm{deg}\;\left(\mathcal{M}^{klmr}_{m^-M_-oM_+m^+p}\Big\vert \mathcal{M}^{klmr}_{m^-M_-oM_+m^+p;1}\right)=\prod_{j=1}^k(m_j^-+1)\prod_{j=1}^m(m_j^++1),
$$
so that the following multiplicative formula holds
$$
\mathrm{deg}\;\left(\mathcal{M}^{klmr}_{m^-M_-oM_+m^+p}\Big\vert \mathcal{M}^{klmr}_{m^-M_-oM_+m^+p;1}\right)\mathrm{deg}\;\mathcal{M}^{klmr}_{m^-M_-oM_+m^+p;1}=\mathrm{deg}\;\mathcal{M}^{klmr}_{m^-M_-oM_+m^+p}.
$$

\medskip

\textbf{Stable region.} We now relate the above computation to stability. The stable region is
$$
\mathcal{E}^n=\bigcup_{l=0}^n\bigcup_{\substack{M_+,M_-,o\\M_++M_-+\vert o \vert=n}}\mathcal{E}^l_{M_-oM_+}
$$
where
$$
\mathcal{E}^l_{M_-oM_+}:=\mathcal{M}^{0l00}_{\mathbf{0}M_-oM_+\mathbf{0}\mathbf{0}},
$$
with $\mathbf{0}=(0,\dots,0)$. This is naturally a stratified space, as follows. Its top open stratum is $\mathcal{E}^n_{0\mathbf{1}0}$ with $\mathbf{1}=(1,\dots,1)$, which is $n$-dimensional, and represented by matrices with only simple elliptic eigenvalues in their spectrum. Its boundary is
$$
\partial \mathcal{E}^n_{0\mathbf{1}0}=\bigcup_{l=0}^{n-1}\bigcup_{\substack{M_+,M_-,o\\M_++M_-+\vert o \vert=n}}\mathcal{E}^l_{M_-oM_+},
$$
represented by matrices having only elliptic eigenvalues, at least one with multiplicity higher than $1$, and/or with $\pm 1$ in the spectrum. The dimension of the boundary stratum $\mathcal{E}^l_{M_-oM_+}$ is 
$$
\dim \;\mathcal{E}^l_{M_-oM_+}=n - M_--M_+-\sum_{j=1}^l(o_j-1)=l.
$$
In particular, $\mathcal{E}^0_{M_-\mathbf{0}M_+}$ is a point (represented by the matrix $\mathds 1_{M_-}\oplus \mathds 1_{M_+}$). Moreover, in order to describe the boundary of $\mathcal{E}^l_{M_-oM_+}$, we introduce a \emph{basic} contraction operator $C_j: \mathbb R^k\rightarrow \mathbb R^{k-1}$ for each $j=1,\dots, l-1$. Given an ordered tuple of real numbers $(a_1,\dots,a_k) \in \mathbb R^{k}$, its contraction on consecutive entries is
$$
C_j(a_1,\dots,a_k):=(a_1,\dots,\wick{\c a_j,\c a_{j+1}}, \dots, a_n):=(a_1,\dots, a_j+a_{j+1},\dots, a_k)\in \mathbb R^{k-1}.
$$
A \emph{contraction} of $a=(a_1,\dots,a_k)$ is by definition the result of applying a sequence of basic contractions to $a$. 

Then the boundary of $\mathcal{E}^l_{M_-oM_+}$ is given by
$$
\partial \mathcal{E}^l_{M_-oM_+}=\bigcup_{l'=0}^{l-1}\bigcup_{\substack{(M'_-,o',M'_+)\in\\ C(M_-,o,M_+)}}\mathcal{E}^{l'}_{M'_-o'M'_+},
$$
where $C(M_-,o,M_+)$ denotes the (finite) set of all possible contractions of $(M_-,o,M_+)\in \mathbb R^{l+2}$.

\medskip

\textbf{The associahedron.} The boundary combinatorics of the stable region can be alternatively encoded as follows. We identify the simple eigenvalues
$$
-1<\nu_1<\dots<\nu_l<1
$$
with the word $-1\nu_1\dots \nu_l 1$. When moving from one stratum to the next, we add a bracket, only allowing to group two consecutive elements at a time. For instance
$$
-1\nu_1\dots \nu_l 1\mapsto -1\nu_1\dots \{\nu_j,\nu_{j+1}\}\dots\nu_l 1
$$
indicates that the eigenvalues $\nu_j$ and $\nu_{j+1}$ come together into a multiplicity two eigenvalue, and so corresponds to the contraction of multiplicities given by
\begin{equation*}
(1,\dots, 1) \mapsto (1,\dots,2,\dots,1).   \end{equation*}
Similarly, a further parenthesis
$$
-1\nu_1\dots \nu_{j-1}\{\nu_j,\nu_{j+1}\}\dots\nu_l 1\mapsto -1\nu_1\dots \{\nu_{j-1},\nu_j,\nu_{j+1}\}\dots\nu_l 1
$$
indicates that the eigenvalue $\nu_{j-1}$ came together with the previous multiplicity two eigenvalue into a multiplicity three eigenvalue, and so corresponds to the contraction
$$
(1,\dots, 1,2,\dots,1)\mapsto (1,\dots,3,\dots,1).
$$
This construction iterates in the obvious way. Here we also allow eigenvalues to come together with $\pm 1$, i.e.\ $\{-1,\nu_1\}\nu_2\{\nu_3,\nu_4,1\}$ is a valid expression. Note that we use the bracket notation to indicate that the order of the elements in the bracket is irrelevant.
Iterating this construction results in a poset of strings (in which all open brackets are accompanied by a corresponding closed one, and there are no nested brackets), and where two strings $a,b$ satisfy $a\leq b$ if $b$ is obtained by a sequence of brackets operations from $a$. This poset then encodes the boundary combinatorics of the stable region, by construction.

\begin{figure}
    \centering
    \includegraphics[width=1\linewidth]{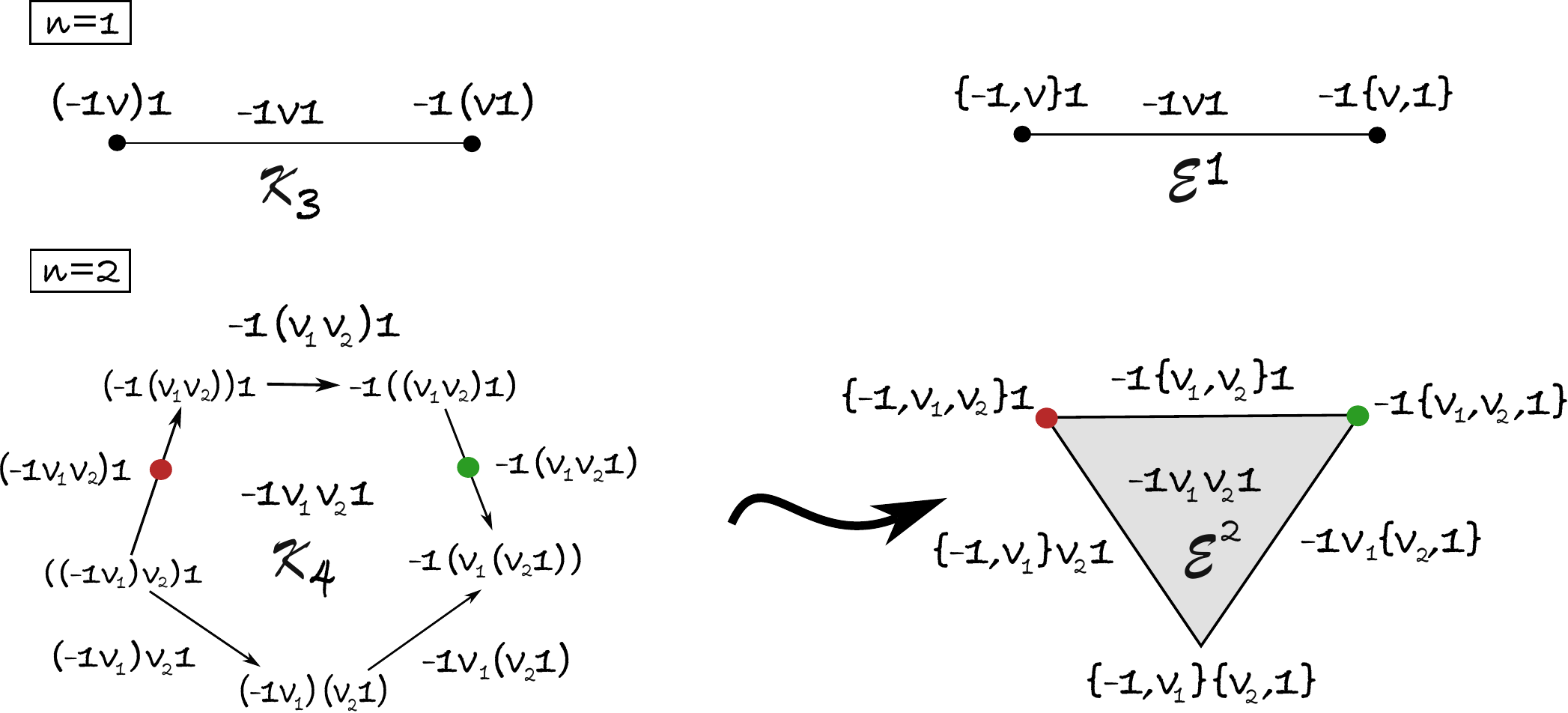}
    \caption{The associahedron $K_3$ coincides with the $1$-dimensional stable region $\mathcal{E}^1$, an interval. The $2$-dimensional stable region (a triangle) is obtained up to homeomorphism by collapsing two edges in $K_4$ (a pentagon) to a point, indicated by points in the Figure.}
    \label{fig:stable_regionlowd}
\end{figure}

Now, the above is intimately related to the operation of taking parenthesis operations
$$
-1\nu_1\dots \nu_l 1\mapsto -1\nu_1\dots (\nu_j\nu_{j+1})\dots\nu_l 1
$$
and iterating them, similarly as above, e.g.\ as
$$
-1\nu_1\dots \nu_{j-1}(\nu_j\nu_{j+1})\dots\nu_l 1\rightarrow -1\nu_1\dots (\nu_{j-1}(\nu_j\nu_{j+1}))\dots\nu_l 1,
$$
and so on, where now a valid expression is for instance $((-1\nu_1)\nu_2)\nu_3(\nu_41)$. The bracket is then the result of removing all interior parenthesis in an expression, i.e.\ symbolically via $(\dots (\dots )\dots)\mapsto (\dots)$, and modding out by the action of the corresponding permutation group (i.e.\ acting on the number of elements inside the bracket), symbolically via
$$
(\underbrace{\dots}_m)\mapsto \{\underbrace{\dots}_m\}=(\underbrace{\dots}_m)/S_m.
$$ For instance, the above expression becomes $\{-1,\nu_1,\nu_2\}\nu_3\{\nu_4,1\}$, where now the order of the elements inside the bracket is irrelevant.

But the combinatorics of expressions with parenthesis is governed by a very-well known polytope, called the \emph{associahedron}. This is the $(m -2)$-dimensional convex polytope $K_m$ in which each vertex corresponds to a way of \emph{correctly} inserting opening and closing parentheses in a string of $m$ letters (meaning that it uniquely determines the order of the product operations), and the edges correspond to single application of the associativity rule. This can also be viewed as a poset, when the arrow indicates that the parentheses have been moved to the right (this is the \emph{Tamari lattice}). Moreover, one can also label the edges with ``\emph{incorrect}'' expressions, containing the common parentheses to each of its boundary vertices. And one can then label the faces by also ``incorrect'' expressions containing the parentheses common to all its boundary edges. And this process iterates, ending in the top strata, which is labelled by the string with no brackets $-1\nu_1\dots\nu_l1$.  

In order to obtain the stable region from the associahedron, we observe that many labels in the latter are actually equivalent when written with the bracket notation. We then conclude:

\medskip

\textbf{Motto.} \emph{The boundary combinatorics of the stable region is determined by the associahedron $K_{n+2}$ in $n+2$ letters, up to collapsing a collection of strata to lower dimensional strata.}

\medskip

\begin{figure}
    \centering
    \includegraphics[width=1.02\linewidth]{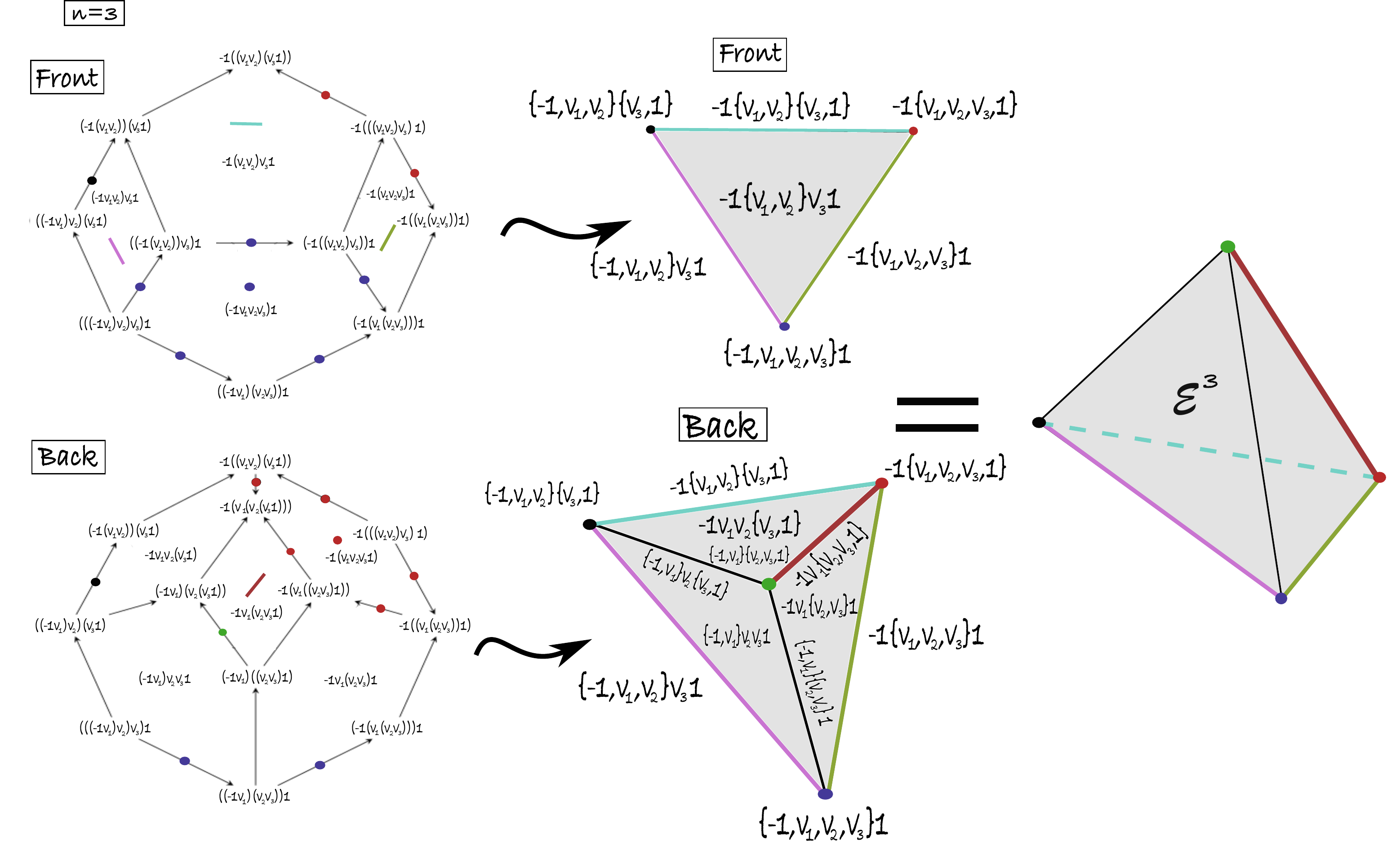}
    \caption{The $3$-dimensional stable region (a tetrahedron) is obtained up to homeomorphism from $K_5$ (an enneahedron) by collapsing 2 faces and 12 edges of $K_5$ to a point --labelled with a point--, and 4 faces of $K_5$ to a line --labelled with a line--.}
    \label{fig:associahedron}
\end{figure}

In other words, the stable region is homeomorphic to a quotient of the associahedron, where we identify those strata whose label become equivalent when written in the bracket notation. The low dimensional cases ($n=1,2,3$) are depicted in Figures \ref{fig:stable_regionlowd} and \ref{fig:associahedron}.

\medskip

\textbf{Branching structure.} We will denote the branches of $\mathcal{E}^l_{M_-oM_+}$, labelled by the $B$-signature $$\epsilon=((c_1,d_1),\dots,(c_l,d_l)),$$ by $\mathcal{E}^{l;\epsilon}_{M_-oM_+},$ so that
$$
\mathrm{int}\;\mathcal{E}^l_{M_-oM_+}=\bigsqcup_\epsilon
\mathcal{E}^{l;\epsilon}_{M_-oM_+}.
$$
By the previous computation of degrees, the interior of the top stratum $\mathcal{E}^n_{0\mathbf{1} 0}$ has $2^n$ branches. Some of these come together at the boundary strata. Which branches can come together over each boundary stratum is determined by the corresponding $B$-signatures, as these behave continuously whenever defined. 

Indeed, to explain this, we first define a few operations. A \emph{middle contraction} operation on consecutive entries of a $B$-signature, symbolically, via the formula
\begin{equation*}
((c_1,d_1),\dots,\wick{(\c c_j,d_j), (c_{j+1},\c d_{j+1})},\dots,(c_l,d_l)):=((c_1,d_1),\dots,(c_j+c_{j+1},d_j+d_{j+1}),\dots,(c_l,d_l)).
\end{equation*}
We also allow the \emph{elimination} contraction of the first and last entries, given by
$$
(\wick{(\c c_1,\c d_1)},\dots,(c_l,d_l))=((c_2,d_2),\dots, (c_{l},d_{l})),
$$
$$
(( c_1, d_1),\dots,\wick{\c( c_l,\c d_l)})=((c_1,d_1),\dots, (c_{l-1},d_{l-1})).
$$
The first operation is the \emph{left} elimination contraction, and the second, the \emph{right} elimination contraction. Then a \emph{contraction} of a $B$-signature is by definition a sequence of middle contractions of consecutive entries, composed with elimination contractions. The \emph{right/left elimination order} of the contraction is the number of right/left elimination contractions in the sequence, the \emph{middle order} is the number of middle contractions in the sequence, and the \emph{total order} is the sum of right/left/middle contractions. Intuitively, a middle contraction represents elliptic eigenvalues of possibly high multiplicity coming together but staying elliptic, a left elimination contraction represents an elliptic eigenvalue of possibly high multiplicity becoming $-1$, and a right elimination contraction represents an elliptic eigenvalue of possibly high multiplicity becoming $1$. In the first case, signatures simply add up, while in the second cases, as the $B$-signature is not defined (the $B$-block becomes degenerate), it is simply suppressed. Note that a contraction of $B$-signatures determines an underlying contraction of multiplicities, by replacing the pair $(c_j,d_j)$ in the above tuple by $o_j=c_j+d_j$. 

In keeping with the bracket notation introduced above, we will sometimes also use the notation $(c,d)=\{\underbrace{+,\dots,+}_{c},\underbrace{-,\dots,-}_{d}\}$, and denote contractions as
\begin{equation*}
\begin{split}
(\{\underbrace{+,\dots,+}_{c_1},\underbrace{-,\dots,-}_{d_1}\},\dots, \wick{\{\underbrace{+,\dots,+}_{c_j}\c,\underbrace{-,\dots,-}_{d_j}\}, \{\underbrace{+,\dots, +}_{c_{j+1}}\c,\underbrace{-,\dots,-}_{d_{j+1}}\}}, \dots, \{\underbrace{+,\dots,+}_{c_l},\underbrace{-,\dots,-}_{d_l}\})
\end{split}
\end{equation*}
$$
:=(\{\underbrace{+,\dots,+}_{c_1},\underbrace{-,\dots,-}_{d_1}\},\dots, \{\underbrace{+,\dots,+}_{c_j+c_{j+1}},\underbrace{-,\dots,-}_{d_j+d_{j+1}}\}, \dots, \{\underbrace{+,\dots,+}_{c_l},\underbrace{-,\dots,-}_{d_l}\}).
$$

\begin{figure}
    \centering
    \includegraphics[width=0.8\linewidth]{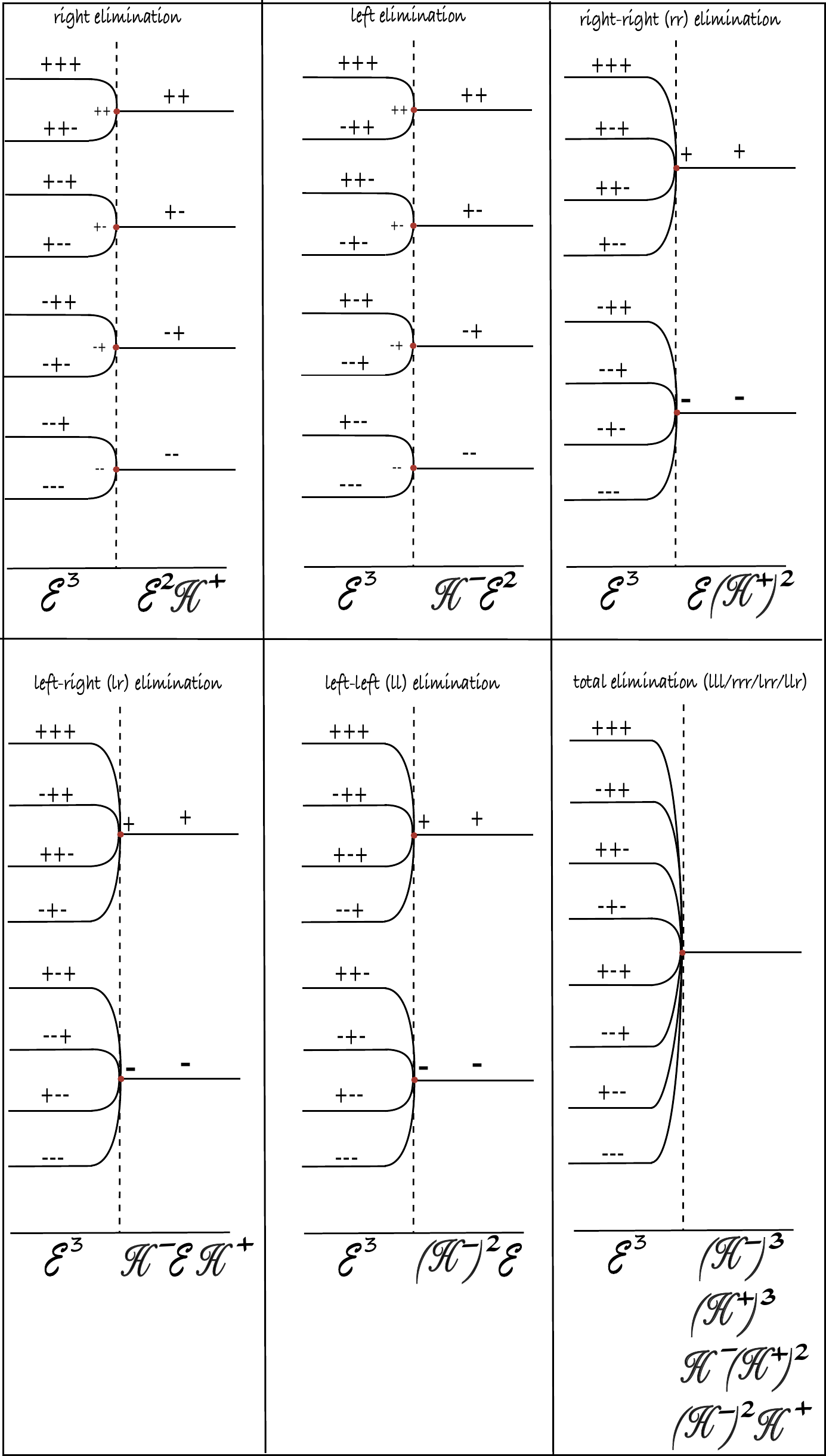}
    \caption{The branching structure for $Sp(6)//Sp(6)$ is obtained from the ones in Figure \ref{fig:eliminations} by collapsing all branches corresponding to hyperbolic eigenvalues together.}
    \label{fig:eliminations_middle}
\end{figure}

Now, the branch $\mathcal{E}^{l;\epsilon}_{M_-oM_+}$ can meet the branch $\mathcal{E}^{l';\epsilon'}_{M_-'o'M_+'}$ at the boundary, where $M'_-o'M'_+=c(M_-,o,M_+)$ for some contraction $c$, if and only if 
$$
\epsilon'=((c'_1,d'_1),\dots, (c'_{l'},d'_{l'}))
$$
is obtained by the underlying contraction of $\epsilon$ determined by $c$, with $$
l'=l-\mathrm{tord}(c),
$$ 
$$
M_-'=M_- + \mathrm{lord}(c),
$$ 
$$
M_+'=M_+ + \mathrm{rord}(c),
$$ 
$$
\vert o'\vert =\vert o \vert-\mathrm{lord}(c)-\mathrm{mord}(c),
$$ 
where $\mathrm{mord}(c),\mathrm{lord}(c),\mathrm{rord}(c), \mathrm{tord}(c)$ are respectively the middle/left/right/total order of $c$. In particular, each $c_j'$ is a sum of consecutive $c_j$'s, and similarly each $d_j'$ is a sum of consecutive $d_j$'s. We then see that two branches can come together at the boundary if and only if their corresponding $B$-signatures have a common contraction, and the branch at which they meet (and therefore the boundary strata at which they meet) is determined by this contraction. The co-dimension of this strata is the total order of the contraction, as per the first of the above formulae. The low-dimensional cases are depicted already in Figures \ref{fig:GIT_sequence} ($n=1$) and \ref{fig:GIT_sequence_3D} ($n=2$). For $n=3$, the branching structure corresponding to transitions from the stable region is depicted in Figures \ref{fig:EN} and \ref{fig:eliminations}. The combinatorics for the branching structure for all other transitions between components is obtained analogously. 

\appendix

\section{Stability, the Krein--Moser theorem, and refinements}\label{app:Krein} Now we explain how the GIT sequence \emph{topologically} encodes the (linear) stability of periodic orbits, and compare it to the basic notions of Krein theory, and the Krein--Moser stability theorem. We will also explain how to obtain Theorem \ref{thm:HN}.

We follow the exposition in Ekeland's book \cite{Eke90} (see also \cite{Ab01}). Consider a linear \emph{symplectic} ODE
$$
\dot x=M(t)x,
$$
where $M(t)=JA(t)$, with $A(t)$ symmetric, and $T$-periodic, i.e.\ $A(t+T)=A(t)$ for all $t$, and $J=\left(\begin{array}{cc}
    0 &  \mathds 1\\
    -\mathds 1 & 0
\end{array}\right)$. The solutions are  $x(t)=R(t)x(0)$, where
$R(t)\in \mathrm{Sp}(2n)$ is symplectic and solves $\dot R(t)=M(t)R(t)$, $R(0)=\mathds 1$. 

\smallskip

\begin{definition}\textbf{(stability)}
The ODE $\dot x=JA(t)x$ is \emph{stable} if all solutions remain bounded for all $t\in \mathbb{R}$. It is \emph{strongly} stable if there exists $\epsilon>0$ such that, if $B(t)$ is symmetric and satisfies $\Vert A(t)- B(t)\Vert <\epsilon$, then the ODE $\dot x=JB(t)x$ is stable. 

A symplectic matrix $R$ is \emph{stable} if all its iterates $R^k$ remain bounded for $k\in \mathbb{Z}$, and it is \emph{strongly} stable if there exists $\epsilon>0$ such that all symplectic matrices $S$ with $\Vert R-S \Vert<\epsilon$ are also stable.    
\end{definition}

One can show that the ODE $\dot x=JA(t)x$ is (strongly) stable if and only if $R(T)$ is (strongly) stable \cite[Section 2, Proposition 3]{Eke90}. Moreover, stability is equivalent to $R(T)$ being diagonalizable (i.e.\ all eigenvalues are semi-simple), with its spectrum lying in the unit circle \cite[Section 1, Proposition 1]{Eke90}. 

Now, consider an elliptic pair $\{\lambda,\overline{\lambda}\}$ of eigenvalues of a symplectic matrix $R$. Then any other symplectic matrix close to $R$ will also have simple eigenvalues in the unit circle different from $\pm 1$ (otherwise an eigenvalue would have to bifurcate into two, as every eigenvalue comes in quadruples, which is not possible if eigenspaces are $1$-dimensional). Therefore in this situation, $R$ is strongly stable. The case of eigenvalues with higher multiplicity is dealt with via Krein theory. Whenever two elliptic eigenvalues come together, this gives a criterion for when they cannot possibly escape the circle and transition into a complex quadruple. This works as follows.

Consider the nondegenerate bilinear form $G(x,y)=x^t \cdot (-iJ) \cdot \overline{y}$ on $\mathbb{C}^{2n}$, associated to the Hermitian matrix $-iJ$. Every real symplectic matrix $R$ preserves $G$. Moreover, if $\lambda,\mu$ are eigenvalues of $R$ which satisfy $\lambda \overline{\mu}\neq 1$, then the corresponding eigenspaces are $G$-orthogonal, since
$$
G(x,y)=G(Rx,Ry)=\lambda \overline{\mu}G(x,y), 
$$
if $x,y$ are the corresponding eigenvectors. Moreover, if we consider the generalized eigenspaces 
$$
E_\lambda=\bigcup_{m\geq 1}\ker(R-\lambda I)^m,
$$
then it also holds that $E_\lambda,E_\mu$ are $G$-orthogonal if $\lambda \overline{\mu}\neq 1$ \cite[Section 2, Proposition 5]{Eke90}. This, in particular, implies that if $\vert\lambda\vert \neq 1$, then $E_\lambda$ is $G$-isotropic, i.e.\ $G\vert_{E_\lambda}=0$. If $\sigma(R)$ denotes the spectrum of $R$, we have a $G$-orthogonal decomposition
$$
\mathbb{C}^{2n}=\bigoplus_{\substack{\lambda \in \sigma(R)\\\vert \lambda \vert\geq 1}}F_\lambda,
$$
where $F_\lambda=E_\lambda$ if $\vert \lambda\vert=1$, and $F_\lambda=E_\lambda \oplus E_{\overline{\lambda}^{-1}}$ if $\vert\lambda\vert >1$. Since $G$ is non-degenerate, and the above splitting is $G$-orthogonal, the restriction $G_\lambda=G\vert_{F_\lambda}$ is also non-degenerate. Note that if $\vert \lambda \vert \neq 1$, with algebraic multiplicity $d$, then the $2d$-dimensional space $F_\lambda$ has $E_\lambda$ as a $d$-dimensional isotropic subspace, and hence the signature of $G_\lambda$ is $(d,d)$. On the other hand, if $\vert\lambda\vert=1$, then the non-degenerate form $G_\lambda$ can have any signature. This justifies the following.

\begin{definition}\textbf{(Krein-positivity/negativity)}
If $\lambda$ is an eigenvalue of the symplectic matrix $R$ with $\vert \lambda \vert=1$, then the signature $(p,q)$ of $G_\lambda$ is called the \emph{Krein-type} or \emph{Krein signature} of $\lambda$. If $q=0$, i.e.\ $G_\lambda$ is positive definite, $\lambda$ is said to be \emph{Krein-positive.} If $p=0$, i.e.\ $G_\lambda$ is negative definite, $\lambda$ is said to be \emph{Krein-negative.} If $\lambda$ is either Krein-negative or Krein-positive, we say that it is \emph{Krein-definite}. Otherwise, we say that it is \emph{Krein-indefinite}.
\end{definition}

If $\lambda$ is of Krein-type $(p,q)$, then $\overline{\lambda}$ is of Krein-type $(q,p)$ \cite[Section 2, Lemma 9]{Eke90}. If $\lambda$ satisfies $\vert \lambda \vert=1$ and it is not semi-simple, then it is easy to show that it is Krein-indefinite \cite[Section 2, Proposition 7]{Eke90}. Moreover, $\pm 1$ are always Krein-indefinite if they are eigenvalues, as they have real eigenvectors $x$, which are therefore $G$-isotropic, i.e.\ $G(x,x)=0$. The following, originally proved by Krein in \cite{Kre1,Kre2,Kre3,Kre4} and independently rediscovered by Moser in \cite{M78}, gives a characterization of strong stability in terms of Krein theory:
\begin{thm}[\textbf{Krein--Moser}]\label{Kreinthm}
 $R$ is strongly stable if and only if it is stable and all its eigenvalues are Krein-definite. 
\end{thm}
See \cite[Section 2, Theorem 3]{Eke90} for a proof. Note that this generalizes the case where all eigenvalues are simple, different from  $\pm 1$ and in the unit circle, as discussed above. Now, the way that the GIT sequence ties up with Krein theory is the following.

\begin{proposition}[\cite{FM}]\label{prop:Krein_vs_B}
    For a Wonenburger matrix, the Krein signature coincides with the $B$-signature, for elliptic eigenvalues.
\end{proposition}

\begin{example}
    As a simple example, to illustrate Proposition \ref{prop:Krein_vs_B}, consider the Wonenburger matrices
$$M=\left(\begin{array}{cccc}
\cos \theta &      0       & -\sin \theta& 0\\
0             & \cos\theta &            0   & -\sin \theta\\
\sin \theta&     0        & \cos \theta  & 0\\
0             & \sin\theta&            0   & \cos \theta
\end{array}\right), 
N=\left(\begin{array}{cccc}
\cos \theta &      0       & \sin \theta & 0\\
0             & \cos\theta &            0   & -\sin \theta\\
-\sin \theta &     0        & \cos \theta & 0\\
0             & \sin\theta &            0   & \cos \theta
\end{array}\right),$$ 
$$
P=\left(\begin{array}{cccc}
\cos \theta&      0       & \sin \theta & 0\\
0             & \cos\theta &            0   & \sin \theta\\
-\sin \theta &     0        & \cos \theta  & 0\\
0             & -\sin\theta &            0   & \cos \theta
\end{array}\right),$$

with $\theta \in (0,\pi)$. These are the three normal forms for $4\times 4$ Wonenburger matrices which are doubly elliptic and have eigenvalues $e^{\pm i \theta}$ of multiplicity two, see \cite{FM}. For the matrix $M$, the eigenspace for the eigenvalue $e^{i \theta}$ is spanned by the eigenvectors
$$
v_1=(1,0,i,0),\; v_2=(0,1,0,i),
$$
and that of $e^{-i \theta}$, by their conjugates 
$$
w_1=(1,0,-i,0),\; w_2=(0,1,0,-i).
$$
With $G=\left(\begin{array}{cc}
   0  &  -i\mathds 1\\
   i\mathds 1  & 0
\end{array} \right)$, we have
$$
G(w_1,w_1)=G(w_2,w_2)=-G(v_1,v_1)=-G(v_2,v_2)=2,
$$
$$
G(w_1,w_2)=G(v_1,v_2)=0.
$$
Therefore the Krein signature of $e^{i \theta}$ is $(0,2)$, and that of $e^{-i\theta}$ is $(2,0)$. By inspection, we see that these are also the $B$-signatures of these eigenvalues. The remaining matrices are dealt with similarly.   
\end{example}

As a corollary of the Krein--Moser theorem and of Proposition \ref{prop:Krein_vs_B}, we obtain the following.

\begin{thm}
 Let $R$ be a Wonenburger matrix. Then $R$ is strongly stable if and only if it is stable and all its eigenvalues are $B$-definite. 
\end{thm}

Having this in mind, the Krein--Moser theorem is detected topologically by the GIT sequence in a very simple and visual way. Indeed, e.g.\ for the case $n=2$, if the stability point lies in the interior of the doubly elliptic region $\mathcal{E}^2$, then the corresponding (equivalence class of) matrices are strongly stable (e.g.\ because they cannot be perturbed away from $\mathcal{E}^2$, but also because of the Krein--Moser theorem). Now, at the boundary of $\mathcal{E}^2$ is where it becomes interesting, and in particular, along the boundary component of $\mathcal{E}^2$ lying in the parabola $\Gamma_d$. By inspecting Figure \ref{fig:GIT_sequence_3D}, we see that the $++$ and the $--$ branches over $\mathcal{E}^2$, in either the middle layer $Sp(4)//Sp(4)$ or the upper layer $Sp^{\mathcal{I}}(2n)/GL_n(\mathbb R)$, do \emph{not} cross from $\mathcal{E}^2$ to $\mathcal{N}$. This means that the boundary of these branches over the portion of $\Gamma_d$ lying over the interval $(-2,2)$ corresponds to two elliptic eigenvalues coming together, and such that the corresponding matrices are \emph{also} strongly stable. This is coherent with the Krein--Moser theorem, as the $B$-signature or equivalentely the Krein signature is positive (respectively negative) definite along the $++$ (respectively the $--$) branch. Note here that we need to use the definition of $B$-signature rather than of $B$-sign, as at $\mathcal{E}^2\cap \Gamma_d$ the eigenvalues are no longer simple.

Moreover, the same phenomenon happens when two positive/negative-hyperbolic eigenvalues of a Wonenburger matrix come together. Indeed, from Figure \ref{fig:GIT_sequence_3D}, we see that one cannot cross from $\mathcal{H}^{\pm}$ to $\mathcal{N}$ if we move along the $++$ or $--$ branch of the top layer. In other words, we obtain the following.

\begin{proposition}
Consider a Wonenburger matrix $M\in Sp^\mathcal{I}(4)$ with a hyperbolic eigenvalue of multiplicity $2$. Then $M$ cannot be perturbed to a Wonenburger matrix with a complex quadruple if and only if its $B$-signature is definite. 
\end{proposition}

In higher dimensions, whether or not a given high-multiplicity elliptic or hyperbolic eigenvalue of a Wonenburger matrix can be perturbed to be a complex quadruple is determined by whether or not its $B$-signature is definite; see e.g.\ Figure \ref{fig:eliminations} and Remark \ref{rk:Krein-Moser}. This gives a topological proof of the Krein--Moser theorem in all dimensions, and in fact generalizes it for the hyperbolic case, in the case of Wonenburger matrices, proving Theorem \ref{thm:HN} in the Introduction.

\end{document}